
\def\krejci{Krej\v c\'\i}

%
%


\def\gianni{\color{red}}
\let\gianni\relax
\def\pier{\color{magenta}}
\let\pier\relax
\def\juerg{\color{green}}
\let\juerg\relax
\def\paolo{\color{blue}}
\let\paolo\relax
\def\ultime{\color{red}}
\let\ultime\relax
\def\corr{\color{red}}
\let\corr\relax

\def\LaTeX{%
  \let\Begin\begin
  \let\End\end
  \def\Bcenter{\Begin{center}}
  \def\Ecenter{\End{center}}
  \let\Label\label
  \let\salta\relax
  \let\finqui\relax
  \let\futuro\relax}

\def\UK{\def\our{our}\let\sz s}
\def\USA{\def\our{or}\let\sz z}



\LaTeX

\USA


\salta

\documentclass[twoside,12pt]{article}
\setlength{\textheight}{24cm}
\setlength{\textwidth}{16cm}
\setlength{\oddsidemargin}{2mm}
\setlength{\evensidemargin}{2mm}
\setlength{\topmargin}{-15mm}
\parskip2mm


\usepackage{color}
\usepackage{amsmath}
\usepackage{amsthm}
\usepackage{amssymb}

\usepackage{amsfonts}
\usepackage{mathrsfs}

\usepackage[mathcal]{euscript}





%

\finqui

\def\Beq{\Begin{equation}}
\def\Eeq{\End{equation}}
\def\Bsist{\Begin{eqnarray}}
\def\Esist{\End{eqnarray}}

\def\Bthm{\Begin{theorem}}
\def\Ethm{\End{theorem}}

\def\Brem{\Begin{remark}\rm}
\def\Erem{\End{remark}}

\def\Bnot{\Begin{notation}\rm}
\def\Enot{\End{notation}}

\let\non\nonumber




\def\step #1 \par{\medskip\noindent{\bf #1.}\quad}


\def\holder{H\"older}
\def\aand{\quad\hbox{and}\quad}

\def\lhs{left-hand side}
\def\rhs{right-hand side}
\def\sfw{straightforward}
\def\wk{well-known}


\def\discretiz{discreti\sz}

\def\organiz{organi\sz}


\def\multibold #1{\def\arg{#1}%
  \ifx\arg\pto \let\next\relax
  \else
  \def\next{\expandafter
    \def\csname #1#1#1\endcsname{{\bf #1}}%
    \multibold}%
  \fi \next}

\def\pto{.}

\def\multical #1{\def\arg{#1}%
  \ifx\arg\pto \let\next\relax
  \else
  \def\next{\expandafter
    \def\csname cal#1\endcsname{{\cal #1}}%
    \multical}%
  \fi \next}


\def\multimathop #1 {\def\arg{#1}%
  \ifx\arg\pto \let\next\relax
  \else
  \def\next{\expandafter
    \def\csname #1\endcsname{\mathop{\rm #1}\nolimits}%
    \multimathop}%
  \fi \next}

\multibold
qwertyuiopasdfghjklzxcvbnmQWERTYUIOPASDFGHJKLZXCVBNM.

\multical
QWERTYUIOPASDFGHJKLZXCVBNM.

\multimathop
dist div dom meas sign supp .


\def\accorpa #1#2{\eqref{#1}--\eqref{#2}}
\def\Accorpa #1#2 #3 {\gdef #1{\eqref{#2}--\eqref{#3}}%
  \wlog{}\wlog{\string #1 -> #2 - #3}\wlog{}}


\def\todx{\scriptstyle\searrow}
\def\tosn{\scriptstyle\nearrow}

\def\somma #1#2#3{\sum_{#1=#2}^{#3}}

\def\graffe #1{\mathopen\{#1\mathclose\}}

\def\<#1>{\mathopen\langle #1\mathclose\rangle}
\def\norma #1{\mathopen \| #1\mathclose \|}

\def\iot {\int_0^t}

\def\iO{\int_\Omega}
\def\intQt{\iot\!\!\iO}

\def\dt{\partial_t}
\def\dn{\partial_\nu}

\def\cpto{\,\cdot\,}

\def\checkmmode #1{\relax\ifmmode\hbox{#1}\else{#1}\fi}
\def\aeO{\checkmmode{a.e.\ in~$\Omega$}}
\def\aeQ{\checkmmode{a.e.\ in~$Q$}}


\def\erre{{\mathbb{R}}}




\def\genspazio #1#2#3#4#5{#1^{#2}(#5,#4;#3)}
\def\spazio #1#2#3{\genspazio {#1}{#2}{#3}T0}

\def\L {\spazio L}
\def\H {\spazio H}
\def\W {\spazio W}

\def\C #1#2{C^{#1}([0,T];#2)}


\def\Lx #1{L^{#1}(\Omega)}
\def\Hx #1{H^{#1}(\Omega)}

\def\Cx #1{C^{#1}(\overline\Omega)}

\def\Luno{\Lx 1}
\def\Ldue{\Lx 2}
\def\Linfty{\Lx\infty}

\def\Huno{\Hx 1}
\def\Hdue{\Hx 2}


\def\LQ #1{L^{#1}(Q)}


\let\theta\vartheta

\let\TeXchi\chi                         
\newbox\chibox
\setbox0 \hbox{\mathsurround0pt $\TeXchi$}
\setbox\chibox \hbox{\raise\dp0 \box 0 }
\def\chi{\copy\chibox}


\def\muz{\mu_0}
\def\rhoz{\rho_0}

\def\rhomin{\rho_*}
\def\rhomax{\rho^*}
\def\rhominbis{\rho_\bullet}
\def\rhomaxbis{\rho^\bullet}

\def\normaV #1{\norma{#1}_V}
\def\normaH #1{\norma{#1}_H}
\def\normaW #1{\norma{#1}_W}

\def\mun{\mu_n}

\def\rhon{\rho_n}
\def\munp{\mu_{n+1}}
\def\munpp{\mu_{n+2}}
\def\rhonp{\rho_{n+1}}
\def\rhonpp{\rho_{n+2}}
\def\rhonm{\rho_{n-1}}
\def\xin{\xi_n}
\def\xinp{\xi_{n+1}}
\def\gn{\gamma_n}
\def\gnp{\gamma_{n+1}}
\def\gz{\gamma_0}
\def\zn{z_n}
\def\znp{z_{n+1}}
\def\znpp{z_{n+2}}
\def\znm{z_{n-1}}
\def\pin{\pi_n}
\def\pinp{\pi_{n+1}}
\def\pim{\pi_m}
\def\piz{\pi_0}

\def\rhoN{\rho_N}
\def\rhoNm{\rho_{N-1}}

\def\mum{\mu_m}
\def\rhom{\rho_m}

\def\rhomp{\rho_{m+1}}
\def\xim{\xi_m}

\def\gm{\gamma_m}

\def\dtau{\delta_\tau}
\def\dzn{\dtau z_n}
\def\dmun{\dtau\mun}
\def\drhon{\dtau\rhon}
\def\drhonp{\dtau\rhonp}
\def\dgn{\dtau\gn}

\def\drhom{\dtau\rhom}

\def\overz{\overline z_\tau}
\def\underz{\underline z_\tau}
\def\hz{\hat z_\tau}
\def\overmu{\overline\mu_\tau}
\def\undermu{\underline\mu_\tau}
\def\hmu{\hat\mu_\tau}
\def\overrho{\overline\rho_\tau}
\def\underrho{\underline\rho_\tau}
\def\hrho{\hat\rho_\tau}
\def\overxi{\overline\xi_\tau}

\def\overg{\overline\gamma_\tau}
\def\underg{\underline\gamma_\tau}
\def\hg{\hat\gamma_\tau}


\let\tau h
\let\hat\widehat

\Begin{document}


\title{{\bf Analysis of a time discretization scheme
  for a nonstandard viscous Cahn-Hilliard
  system}\footnote{{\bf Acknowledgments.}\quad\rm
The present paper benefits from the support of the MIUR-PRIN Grant 2010A2TFX2 ``Calculus of variations'' for PC and~GG, the GA\v CR Grant P201/10/2315 
and RVO:~67985840 for~PK, the DFG Research Center {\sc Matheon} 
in Berlin for JS.}}
\author{}
\date{}
\maketitle

\Bcenter
\vskip-2.2cm
{\large\bf Pierluigi Colli$^{(1)}$}\\
{\normalsize e-mail: {\tt pierluigi.colli@unipv.it}}\\[.4cm]
{\large\bf Gianni Gilardi$^{(1)}$}\\
{\normalsize e-mail: {\tt gianni.gilardi@unipv.it}}\\[.4cm]
{\large\bf Pavel \krejci $^{(2)}$}\\
{\normalsize e-mail: {\tt krejci@math.cas.cz}}\\[.25cm]
{\large\bf Paolo Podio-Guidugli$^{(3)}$
}\\
{\normalsize e-mail: {\tt ppg@uniroma2.it}}\\[.25cm]
{\large\bf J\"urgen Sprekels$^{(4)}$}\\
{\normalsize e-mail: {\tt sprekels@wias-berlin.de}}\\[.6cm]
$^{(1)}$
{\small Dipartimento di Matematica ``F. Casorati'', Universit\`a di Pavia}\\
{\small Via Ferrata 1, 27100 Pavia, Italy}\\[.2cm]
$^{(2)}$
{\small Institute of Mathematics, Academy of Sciences of the Czech Republic}\\
{\small Zitna 25, CZ - 115 67 Praha 1, Czech Republic}\\[.2cm]
$^{(3)}$
{\small Dipartimento di Ingegneria Civile e Ingegneria Informatica}\\
{\small Universit\`a di Roma {\ultime TorVergata},
Via Politecnico 1, 00133 Roma, Italy}\\[.2cm]
$^{(4)}$
{\small Weierstra\ss-Institut f\"ur Angewandte Analysis und Stochastik}\\
{\small Mohrenstra\ss e\ 39, 10117 Berlin, Germany}\\[.8cm]
\Ecenter

{
\Begin{abstract}
{\paolo In t}his paper {\paolo we propose a} time discretization 
of a system of two parabolic equations {\paolo describing diffusion-driven atom rearrangement in crystalline matter. The equations express the balances
of microforces and microenergy; the two phase fields are the order parameter
and the chemical potential.} The initial and boundary{\paolo -}value problem 
for the evolutionary system is known to be well posed. {\paolo Convergence of the discrete scheme to the solution of 
the continuous problem is proved} by a careful 
development of uniform estimates, by weak compactness and a suitable treatment of
nonlinearities. {\pier Moreover\paolo , for} the difference of discrete
and continuous solutions {\paolo we} prove an error estimate of order one with respect to the time step. 
\\[2mm]
\\[2mm]
\End{abstract}
}
\salta

\pagestyle{myheadings}
\newcommand\testopari{\sc Colli \ --- \ Gilardi \ --- \ \krejci \ --- \ Podio-Guidugli \ --- \ Sprekels}
\newcommand\testodispari{\sc Time discretization of a nonstandard viscous Cahn-Hilliard system}
\markboth{\testodispari}{\testopari}


\finqui


\section{Introduction}
\label{Intro}
\setcounter{equation}{0}

This paper deals with a time discretization of a PDE system arising from a mechanical model for phase segregation {\paolo by} atom rearrangement on a lattice. 
The model was proposed by one of us in \cite{Podio}; the resulting system has 
been analyzed in a recent and intensive research work by four of the {\paolo present} authors {\paolo (see, in particular,  \cite{CGPS3} and \cite{CGPS6}, both for well-posedness results and for a detailed presentation of the model).} 

The initial and boundary value problem we consider consists in looking 
for two \emph{fields},  the {\it chemical potential} $\mu>0$ and the \emph{order 
parameter} $\rho\in (0,1)$, solving 
\Bsist
  & (1 + 2g(\rho)) \dt\mu + \mu \, \dt g(\rho) - \Delta\mu = 0
  & \quad \hbox{in $\Omega \times (0,T),$}
  \label{Iprima}
  \\
  & \dt\rho - \Delta\rho + f'(\rho) = \mu g'(\rho)
  & \quad \hbox{in $\Omega \times (0,T),$}
  \label{Iseconda}
  \\
  & \dn\mu = \dn\rho = 0
  & \quad \hbox{on $\Gamma \times (0,T),$}
  \label{Ibc}
  \\
  & \mu|_{t=0} = \muz
  \aand
  \rho|_{t=0} = \rhoz
  & \quad \hbox{in $\Omega$};
  \label{Icauchy}
\Esist
\Accorpa\Ipbl Iprima Icauchy
here $\Omega$ denotes a bounded domain of $ \erre^3$ with conveniently smooth boundary $\Gamma$, $T>0$ stands for some final time, and $\dn$ denotes differentiation in the direction of the outward normal $\nu$. 

{\paolo Problem \Ipbl\ } is parameterized by two nonlinear scalar-valued functions, $g$ and $f$, which enter into the definition of the system's \emph{free energy}:
\begin{equation}
\label{fe-2}
\psi=\widehat\psi(\rho,\nabla\rho,\mu) 
= - {}{\pier \left(\frac1 2 + g(\rho)\! \right)\! \mu}  + f (\rho) +\frac{1}{2}|\nabla\rho|^2.
\end{equation}
We point out that in {\paolo \eqref{Iprima}--\eqref{fe-2}} {\juerg all} physical constants have been set equal to $1${\paolo . We also note that the last two terms in \eqref{fe-2}} favor phase segregation, the former because it introduces local energy minima, the latter because it penalizes spatial changes of the order parameter. For $g$, one can take any smooth function, provided it is nonnegative in the {\paolo physically admissible domain}:
\begin{equation}\label{newh}
g(\rho) \geq 0 \quad {\pier \hbox{for all  $\, \rho\in (0,1)$;}}
\end{equation}
{\paolo accordingly,} the coefficient $1/2$ of $\mu$ in \eqref{fe-2} {\paolo should be regarded as a prescribed material} bound.  
%
%
{\paolo As to} the possibly multi-well potential $f$, we {\paolo take it to be} the sum
of two functions:
\[
f(\rho)=f_1(\rho) + f_2(\rho); 
\] 
{\paolo the one, $f_1$, is convex over $(0,1)$}, and such that its derivative $f'_1$ (and possibly {\juerg also} $f_1$) is singular at the 
endpoints $0$ and $1$ (cf.~{\paolo \eqref{hpfp}}); the other {\paolo is required to be smooth over the entire interval $ [0,1]$, but not to have any convexity property,} so that in 
equation~\eqref{Iseconda}  $f'_2$
may {\paolo serve} as a non-monotone perturbation of the increasing function {\paolo $f_1'$}.

{\paolo As to the parameter functions, in \cite{CGPS3} the choice made for $g$ was:}
\begin{equation}\label{oldh}
g(\rho)=\rho,
\end{equation}
{\paolo while the assumptions on $f$ were compatible with choosing a double-well potential:} 
\begin{equation}
\label{ex1} 
f(\rho)= \alpha_1 \, \{ \rho\,\ln (\rho)+(1-\rho)\,\ln (1-\rho) \} + \alpha_2 \, \rho \, (1-\rho) + \alpha_3 \rho,
\end{equation}
for some non-negative constants $\alpha_1, \,  \alpha_2 ,\, \alpha_3 $. Note that, if {\paolo $\alpha_3$ is taken null, then,} according to whether or not $2\alpha_1 \geq \alpha_2$, $f$ {\paolo turns out to be} convex in the whole of $[0,1]$ or it exhibits two wells, with a local maximum at $\rho=1/2$; {\paolo moreover, for $ \alpha_3>0 $, the combined function:} 
\begin{equation}
\non
- {\pier g(\rho)\mu} + f (\rho)   \quad \hbox{(a part of }\, \psi) 
\end{equation}
shows one global minimum in all cases, and {\juerg it depends on the sign of $(\alpha_3 - \mu)$ which minimum actually occurs}. 
On the other hand, the framework of paper \cite{CGPS6} {\paolo allows for much more general choices of} $g$ and $f $, {\paolo as well as for nonlinear diffusion of} $\mu$. Existence and uniqueness results 
{\juerg were} proved in both \cite{CGPS3}
and \cite{CGPS6}, with different approaches. Here, we take inspiration from arguments developed either in {\paolo the one or in the other of those papers}.    

We introduce a time discretization of system \Ipbl\ which is implicit with respect to the principal terms and tries to handle very carefully the nonlinearities. Namely,
we address the recursive sequence of the elliptic problems:
\Bsist
  & (1+2\gn) \, {\ultime \dtau \mu_n} + \munp \, {\ultime \dtau \gamma_n} - \Delta\munp
  = 0 & \quad \hbox{in $\Omega$,}
  \label{Idprima}
  \\
  & {\ultime \dtau \rho_n} - \Delta\rhonp + f'(\rhonp) 
  = \mun g'(\rhon) & \quad \hbox{in $\Omega $,}
  \label{Idseconda}
  \\
 & 
  \dn\munp = \dn\rhonp = 0
   & \quad \hbox{on $\Gamma $,}
  \label{Idbc}
\Esist
\Accorpa\Idpbl Idprima Idbc
{\pier for $n=0,1, \ldots , N $, where $\tau =T/N$ is the time step,} $\gamma_n:=g(\rho_n)$ and{\ultime , for any $(N+1)$-ple $z_0, z_1, \ldots , z_N $, we let} 
$${\ultime \dtau z_n}:=(z_{n+1}-z_n)/\tau  {\pier \quad\hbox{ for } \,  n= 0, \ldots , N-1.} $$ 
After showing the existence of a discrete solution at any step, we carry out a number of uniform estimates on the time-discrete solution which allow us to prove convergence to the unique solution $(\mu,\rho) $ of the continuous problem \Ipbl, as {\pier $\tau$ tends to $0$ (or, equivalently{\ultime ,} $N$ goes to 
$+ \infty$)}. Then, we estimate {\paolo certain} norms of the difference between the {\paolo piecewise-linear-in-time} interpolants of the discrete solutions and the continuous solution: {\paolo more precisely}, the first {\paolo error estimate we prove is} of order $\tau^{1/2}${\pier ;} the second, {\paolo which holds} under stronger regularity assumptions on the initial data, {\paolo is} of order $\tau.$  

We {\paolo regard} our results {\paolo as a} cornerstone {\paolo in the construction of} a  time-and-space discretization 
of problem {\Ipbl}. {\paolo With reference to such a complete} discretization of Cahn-Hilliard 
and viscous Cahn-Hilliard systems, {\paolo we} quote papers 
\cite{BEGSS95, BB95, BB01, BBG99, BBG01, BM10, BM11, BE92, 
GK09-1, GK09-2, GKS13}. Some recent efforts can be found in the  
literature with the aim of analyzing other classes of phase transition problems, 
either to show existence via time discretization~\cite{B04, C11, CFK01, FKS11, 
GS06, KS10, R05, S06, V07} or to prove numerical results 
such as special convergence properties, stability or error 
estimates~\cite{CP01, CNS00, CPM10, EJK10, KS02, MPP09, S04, S99, S00}~{\juerg (cf. also \cite{K12} 
for a recent review on phase-field models)}. 
{\paolo We dare} say that our contribution goes deeply into the structure of the mathematical problem, 
{\paolo because, as is not the case for many other similar investigations,} we succeed in showing a linear order of convergence.
 
{\paolo Our} paper is \organiz ed as follows.
In the next section, we {\paolo list and discuss our} assumptions, {\paolo formulate the continuous and discrete problems precisely,} and state our main results. Section~\ref{Existence} is devoted to {\paolo proving that there is}  a discrete solution. The convergence result is proved in the long and articulate Section~\ref{Convergenceproof}. Finally, the last two {\paolo Sections~\ref{ProofErrore} and~\ref{ProofErrorebis}} contain detailed proofs of the two error estimates.

%
%


\section{Main results}
\label{MainResults}
\setcounter{equation}{0}

In this section, we describe the mathematical problem under investigation,
introduce the time \discretiz ation scheme, 
make our assumptions precise, and state our results.

First of all,
we assume $\Omega$ to be a bounded connected open set in $\erre^3$
with smooth boundary~$\Gamma$.
{\paolo For convenience,} we set
\Beq
  V := \Huno,
  \quad H := \Ldue ,
  \aand
  W := 
  \graffe{v\in\Hdue:\ \dn v = 0 ,\ \hbox{on $\Gamma$}},
  \label{defspazi}
\Eeq
and we endow these spaces with their standard norms,
for which we use a self-explanato\-ry notation 
like $\normaV\cpto$.
The~notation $\norma\cpto_p$ ($1\leq p\leq+\infty$) 
stands for the standard $L^p$-norm in~$\Lx p$;
{\juerg for short, we sometimes do not distinguish 
between a space (or~its norm) and a power thereof}.

{\paolo As to the parameter functions~$f$ and $g$}, we assume~that
\Bsist
  \hskip-1cm && f = f_1 + f_2,
  \quad \hbox{where}
  \label{hpf}
  \\
  \hskip-1cm && \hbox{$f_1:(0,1)\to[0,+\infty)$ is a convex $C^2$ function satisfying}
  \non 
  \\
  \hskip-1cm && \hskip2.8cm \lim_{r\todx 0} f_1'(r) = - \infty\,,
  \aand
  \lim_{r\tosn 1} f_1'(r) = + \infty ,
  \label{hpfp}
  \\
  \hskip-1cm && f_2:[0,1] \to \erre \quad \hbox{is of class $C^2$;}
  \label{hpfd}
  \\
  \hskip-1cm && g:[0,1] \to \erre \quad \hbox{is of class $C^2$ and nonnegative}.
  \label{hpg}
\Esist
\Accorpa\Hpf hpf hpg
For the initial data, we require that
\Bsist
  && \muz \in V \cap \Linfty
  \aand \muz \geq 0 \quad \aeO;
  \label{hpmuz}
  \\
  && {\ultime \rhoz \in W \subset \Cx0 
  \aand \inf\rhoz > 0 , \ \, \sup\rhoz <1 .}
  \label{hprhoz}
\Esist
\Accorpa\Hpdati hpmuz hprhoz
\Accorpa\Hpfdati hpf hprhoz
{\ultime We stress that conditions \eqref{hprhoz} actually imply that 
$\rhoz$ is $1/2$-H\"older continuous: indeed, as $\Omega$ is a 
three-dimensional domain, $W$ is continuosly embedded in $C^{0,1/2} 
(\/{\overline\Omega}\/)$. As a consequence, also $f(\rhoz)$ and $f'(\rhoz)$ 
are $1/2$-H\"older continuous, since $f$ and $f'$ are smooth in~$(0,1)$. On 
the other hand, we point out that in the sequel we will mostly exploit the 
compactness of the embedding $W \subset  \Cx0 $; H\"older continuity 
will play no role.}

{\paolo As recalled in the Introduction, in papers \cite{CGPS3} and \cite{CGPS6} two versions of problem~\Ipbl\ were} solved {\paolo
over an arbitrary time interval~$[0,T]$ in a rather strong sense, because} the solution {\paolo pairs} $(\mu,\rho)$ {\paolo were} required to satisfy
\Bsist
  && \mu \in \H1H \cap \L2W ,
  \label{regmu}
  \\
  && \rho \in \W{1,\infty}H \cap \H1V \cap \L\infty W ,
  \label{regrho}
  \\
  && \mu \geq 0 \quad \aeQ,
  \label{mupos}
  \\
  && 0 < \rho <1 \quad \aeQ
  \aand
  f'(\rho) \in \L\infty H.
  \label{regfprho}
\Esist
\Accorpa\Regsoluz regmu regfprho
Note that the boundary conditions~\eqref{Ibc} follow from
\accorpa{regmu}{regrho}, due to the definition of~$W$ in~\eqref{defspazi}.
{\paolo Accordingly, t}he solution{\paolo s} to the  problem{\paolo s of type  \Ipbl\ studied in \cite{CGPS3} and \cite{CGPS6} were pairs} $(\mu,\rho)$ satisfying, {\paolo in addition to~\Regsoluz,} the system
\Bsist
  & \bigl( 1 + 2g(\rho) \bigr) \dt\mu + \mu \, \dt g(\rho) - \Delta\mu = 0
  & \quad \aeQ,
  \label{prima}
  \\
  & \dt\rho - \Delta\rho + f'(\rho) = \mu g'(\rho)
  & \quad \aeQ,
  \label{seconda}
  \\
  & \mu(0) = \muz
  \aand
  \rho(0) = \rhoz
  & \quad \aeO .
  \label{cauchy}
\Esist
\Accorpa\Pbl prima cauchy
{\paolo Some of the results proved in the quoted papers are summarized in the following theorem.}

\Bthm
\label{Esistcont}
{\paolo Let} {\juerg assumptions \Hpfdati\ hold}.
Then, there exists a unique pair $(\mu,\rho)$ {\ultime satisfying} \Regsoluz\ and solving
problem~\Pbl.
Moreover, $\mu\in\LQ\infty$, and there exist $\rhomin,\rhomax\in(0,1)$ such that
$\rhomin\leq\rho\leq\rhomax$ \aeQ.
\Ethm

The {\paolo main} aim of the present paper is to 
{\paolo show that, given the time-\discretiz ation scheme
introduced here below,}
the discrete solution converges to the solution $(\mu,\rho)$ 
as the time step $\tau$ tends to zero.

\Bnot
\Label{Gennot}
Assume that $N$ is a positive integer, and {\paolo let $Z$ be} any normed space.
We define $\dtau:Z^{N+1}\to Z^N$ as follows:
\Bsist
  && \hbox{for $z=(z_0,z_1,\dots ,z_N)\in Z^{N+1}$ and $w=(w_0,\dots,w_{N-1})\in Z^N$,} 
  \non
  \\
  && \dzn = w
  \quad \hbox{means that} \quad
  w_n {\paolo =} \frac{z_{n+1}-z_n} \tau
  \quad \hbox{for $n=0,\dots,N-1$}.
  \qquad
  \label{defdelta}
\Esist
We can also iterate {\paolo the discretization procedure,
and define, e.g.,} 
\Beq
  \dtau^2\zn := \frac {\dtau\znp - \dtau\zn} \tau
  = \frac {\znpp - 2 \znp + \zn} {\tau^2}
  \quad \hbox{for $n=0,\dots,N-2$}.
  \label{defdeltadue}
\Eeq
Next, by setting $\tau:=T/N$
(without stressing the {\juerg dependence} of $\tau$ on~$N$)
and $I_n:=((n-1)\tau,n\tau)$ for $n=1,\dots,N$,
we introduce the interpolation maps
from $Z^{N+1}$
into either $\L\infty Z$ or $\W{1,\infty}Z$ as follows:
for $z=(z_0,z_1,\dots ,z_N)\in Z^{N+1}$, we~set
\Bsist
  && \overz ,\, \underz \in \L\infty Z 
  \aand 
  \hz \in \W{1,\infty}Z \, ,
  \label{reginterp}
  \\
  && \overz(t) = \zn
  \aand
  \underz(t) = \znm
  \quad \hbox{for a.a.\ $t\in I_n$, \ $n=1,\dots,N$},
  \label{pwconstant}
  \\
  && \hz(0) = z_0
  \aand
  \dt\hz(t) = \dtau\znm
  \quad \hbox{for a.a.\ $t\in I_n$, \ $n=1,\dots,N$}.
  \qquad
  \label{pwlinear}
\Esist
\Enot
%
{\paolo \noindent These} maps yield the {\juerg backward/forward} piecewise{\paolo -}constant
and piecewise{\paolo -}linear interpolants of the discrete vectors.
We obviously have:
\Bsist
  && \norma\overz_{\L\infty Z}
  = \max_{n=1,\dots,N} \norma\zn_Z \,, \quad
  \norma\underz_{\L\infty Z}
   = \max_{n=0,\dots,N-1} \norma\zn_Z,
  \label{ouLinftyZ}
  \\
  && \norma\overz_{\L2Z}^2
  = \tau \somma n1N \norma\zn_Z^2 \,, \quad
  \norma\underz_{\L2Z}^2
  = \tau \somma n0{N-1} \norma\zn_Z^2 \,.
  \label{ouLdueZ}
\Esist
Moreover, as $\hz(t)$ is a convex combination
of $\znm$ and $\zn$ for $t\in I_n$, we also have
\Bsist
  && \norma\hz_{\L\infty Z}
  = \max_{n=1,\dots,N} \max\{\norma\znm_Z,\norma\zn_Z\}
  = \max\{\norma{z_0}_Z,\norma\overz_{\L\infty Z}\},
  \qquad\qquad
  \label{normahzLinftyZ}
  \\
  && \norma\hz_{\L2Z}^2
  \leq \tau \somma n1N \bigl( \norma\znm_Z^2 + \norma\zn_Z^2 \bigr)
  \leq \tau \norma{z_0}_Z^2
  + 2 \norma\overz_{\L2Z}^2 \,.
  \label{normahzLdueZ}
\Esist
Finally, by a direct computation,
it is \sfw\ to prove that
\Bsist
  && \norma{\overz-\hz}_{\L\infty Z}
  = \max_{n=0,\dots,N-1} \norma{\znp-\zn}_Z
  = \tau \, \norma{\dt\hz}_{\L\infty Z},
  \qquad
  \label{interpLinfty}
  \\
  && \norma{\overz-\hz}_{\L2Z}^2
  = \frac \tau 3 \somma n0{N-1} \norma{\znp-\zn}_Z^2
  = \frac {\tau^2} 3 \, \norma{\dt\hz}_{\L2Z}^2,
  \label{interpLdue}
\Esist
and that the same identities hold for the difference $\underz-\hz$.

At this point, we can write the discrete scheme presented in the Introduction
in a precise form.
For any positive integer~$N$, we look for two vectors 
$(\mun)_{n=0}^N$ and $(\rhon)_{n=0}^N$ satisfying the {\paolo following} conditions:

\smallskip
\sl

\indent\llap{$i)$ }
the first components $\muz$ and $\rhoz$ coincide with the initial data;

\smallskip

\indent\llap{$ii)$}
for $n=0,\dots,N-1,$ we have that
\Bsist
  \munp \,,\, \rhonp \in W , \quad
  \munp \geq 0
  \aand
  0 < \rhonp < 1
  \quad \hbox{in $\Omega$}, \quad
  f'(\rhonp) \in H;
  \quad
  \label{regdsoluz}
\Esist

\indent\llap{$iii)$}
{\pier if} $(\gn)_{n=0}^N$ {\pier is the vector} {\paolo whose components are} 
$\gn:=g(\rhon)$, {\pier there hold}
\Bsist
  && (1+2\gn) \, \dmun + \munp \, \dgn - \Delta\munp
  = 0,
  \label{dprima}
  \\
  && \drhon - \Delta\rhonp + f'(\rhonp) 
  = \mun g'(\rhon),
  \label{dseconda}
\Esist
\Accorpa\Discreto dprima dseconda
\indent
for $n=0,\dots,N-1$.

\smallskip
\rm

\noindent Also in this case, the homogenous Neumann boundary conditions
are {\paolo implicit} in the regularity requirements
(see \eqref{regdsoluz} and~\eqref{defspazi}).

Clearly, the ``true'' problem consists in finding $(\munp,\rhonp)$
once $(\rhon,\mun)$ is given.
Here is our result in this direction.

\Bthm
\label{Esistdiscr}
Assume \Hpfdati.
Then, there exists $\tau_0>0$ such that,
for $\tau<\tau_0$ and $n=0,\dots,N-1$, problem~\Discreto\ has
a unique solution $(\munp,\rhonp)$ satisfing \eqref{regdsoluz}.
\Ethm

{\paolo Our} next results {\paolo concern} firstly convergence of interpolants
for vectors $(\rhon)$ and $(\mun)$ to the solution 
$(\mu,\rho)$ to problem \Pbl, then error estimates.
{\paolo We point out that, for simplicity, the convergence theorem here below is not stated in a precise form: the topological setting will be specified later, by means of relations} \futuro\eqref{convmu}-\futuro\eqref{convdtrho}.

\Bthm
\label{Convergenza}
Assume \Hpfdati.
Then,  {\paolo in accord {\pier with}} Notation~\ref{Gennot},
the sequences of interpolants for the discrete solutions
given by Theorem~\ref{Esistdiscr} converge to the solution $(\mu,\rho)$
given by Theorem~\ref{Esistcont} as $\tau$ tends to~$0$, in a {\paolo suitable} topology.
\Ethm

\Bthm
\label{Errore}
In addition to \Hpfdati,
assume that
\Beq
  \rhoz \in \Hx3 .
  \label{hperrore}
\Eeq
Then, {\juerg for sufficiently small $h>0$}, {\paolo the following error estimate holds:}
\Beq
  \norma{\hrho-\rho}_{\H1H\cap\L\infty V}
  + \norma{\hmu-\mu}_{\L\infty H\cap\L2V}
  \leq c \, \tau^{1/2},
  \label{errore}
\Eeq
where $c$ depends {\paolo only} on the structural assumptions and the data.
\Ethm

\Bthm
\label{Errorebis}
{\pier In addition to \Hpfdati}, assume 
{\gianni\eqref{hperrore}} and
\Beq
  \muz \in W .
  \label{hperrorebis}
\Eeq
Then, {\juerg for sufficiently small $h>0$}, {\paolo the following error estimate holds:}
\Beq
  \norma{\hrho-\rho}_{\H1H\cap\L\infty V}
  + \norma{\hmu-\mu}_{\L\infty H\cap\L2V}
  \leq c \, \tau,
  \label{errorebis}
\Eeq
where $c$ depends {\paolo only} on the structural assumptions and the data.
\Ethm

\Brem
\label{Remhperrore}
It is easy to see that our assumptions \Hpfdati\
ensure that both $f'(\rhoz)$ and $\muz g'(\rhoz)$ belong to~$V$.
It follows that \eqref{hperrore} is equivalent~to
\Beq
  -\Delta\rhoz + f'(\rhoz) - \muz g'(\rhoz) \in V .
  \label{dahperrore}
\Eeq
We also notice that {\juerg the} assumptions \eqref{hperrore} and \eqref{hperrorebis}
ensure further regularity for the solution $(\mu,\rho)$ to the continuous problem
(see the forthcoming Remark~\futuro\ref{Piureg}).
\Erem

{\paolo We prove Theorem~\ref{Esistdiscr} in  Section~\ref{Existence} and Theorem~\ref{Convergenza}
in  Section~\ref{Convergenceproof}; the}
last two sections are devoted to {\paolo proving, respectively,} Theorems~\ref{Errore} and~\ref{Errorebis}.

Throughout the paper,
we account for the \wk\ {\juerg embeddings} $V\subset\Lx q$
$(1\leq q\leq 6)$ and $W\subset\Cx0$,
and for the related Sobolev inequalities:
\Beq
  \norma v_q \leq C \normaV v
  \aand
  \norma v_\infty \leq C \norma v_W,
  \label{sobolev}
\Eeq
for $v\in V$ 
and $v\in W$, respectively,
where $C$ depends on~$\Omega$ only,
since sharpness is not needed.
We remark that {\paolo these embeddings 
are compact.}
In particular,  the {\paolo following compactness inequality holds:}
\Beq
  \norma v_4 \leq \sigma \normaH{\nabla v} + C_\sigma \normaH v,
  \quad \hbox{for every $v\in V$ and $\sigma>0$};
  \label{compact}
\Eeq
in \eqref{compact}, $C_\sigma$ is a constant that depends only on $\Omega$ and~$\sigma$.
Furthermore, we {\juerg make repeated} use of
 {\paolo \holder\ inequality,  of 
the following elementary identity: }
\Bsist
  && (a-b) a = \frac 12 \, a^2 - \frac 12 b^2 + \frac 12 \, (a-b)^2,
  \quad \hbox{for every $a,b\in\erre$,}
  \label{elementare}
\Esist
{\paolo and of Young's} inequality
\Bsist
  && ab \leq \sigma a^2 + \frac 1{4\sigma} \, b^2,
  \quad \hbox{for every $a,b\geq 0$ and $\sigma>0$}.
  \label{young}
\Esist
{\juerg Moreover, we use the discrete Gronwall lemma in the following form
(see, e.g., \cite[Prop. 2.2.1]{Jerome}): if
 $(a_0,\dots,a_N)\in[0,+\infty)^{N+1}$ and $(b_1,\dots,b_N)\in[0,+\infty)^N$ satisfy} 
\Bsist
  && a_m \leq a_0 + \somma n1{m-1} b_n a_n 
  \quad \hbox{for $m=1,\dots,N$},
  \quad \hbox{then }
  \non
  \\
  && a_m \leq a_0 \exp \Bigl( \somma n1{m-1} b_n \Bigr)
  \quad \hbox{for $m=1,\dots,N$}.
  \label{dgronwall}
\Esist

Finally, throughout the paper
we use a small-case italic $c$ for {\paolo a number of} different constants that
may only depend on~$\Omega$, the final time~$T$, the shape of~$f$, 
the properties of the data involved in the statements at hand;
those constants we need to refer to are always denoted by 
capital letters, just like $C$ in~\eqref{sobolev}.
Moreover, a~notation like~$c_\sigma$ 
signals a constant that depends also on the parameter~$\sigma$.
The reader should keep in mind that the meaning of $c$ and $c_\sigma$ might
change from line to line and even in the same chain of inequalities
and that their values {\paolo never} depend on the time step~$\tau$.


\section{Existence}
\label{Existence}
\setcounter{equation}{0}

In this section, we prove Theorem~\ref{Esistdiscr}.
We argue inductively with respect to~$n$, i.e.,
by assuming that {\paolo a pair $(\mun,\rhon)$ satisfying~\eqref{regdsoluz}  
with $n$ in place of~$n+1$ is given,}
we prove that problem \Discreto\ has a unique solution $(\munp,\rhonp)$
satisfying~\eqref{regdsoluz}.
More precisely, as {\paolo is going to be} clear from the proof,
we need less regularity for $\mun$,
e.g., $\mun\in V$.
In particular, our assumptions on $\muz$ are sufficient to start.
We rewrite \Discreto\ in the form
\Bsist
  && \bigl( 1+\gn+\gnp \bigr) \munp - \tau \Delta\munp
  = \bigl( 1+2\gn \bigr) \mun,
  \label{dprimabis}
  \\
  && \rhonp - \tau \Delta\rhonp + \tau f'(\rhonp) 
  = \rhon + \tau \mun g'(\rhon),
  \label{dsecondabis}
\Esist
and solve \eqref{dsecondabis} {\paolo first for $\rhonp$ (so~that $\gnp$ is also known),} then~\eqref{dprimabis}.
In~order to solve both problems, {\paolo it is expedient to replace each equation by} a minimum problem,
at least for $\tau$ small enough.
We consider the functionals:
\Bsist
  && J_1 : V \to \erre
  \aand
  J_2 : V \to (-\infty, +\infty],
  \quad \hbox{defined by, {\paolo respectively,}}
  \non
  \\ 
  && J_1(v) 
  := \frac \tau 2 \iO |\nabla v|^2
  + \frac 12 \iO \bigl( 1+\gn+\gnp \bigr) v^2
  - \iO \bigl( 1+2\gn \bigr) \mun \, v{\paolo \aand}
  \qquad
  \label{minmu}
  \\
  && J_2(v) 
  := \frac \tau 2 \iO |\nabla v|^2
  + \frac 12 \iO v^2
  + \tau \iO \tilde f(v)
  - \iO \bigl( \rhon + \tau \mun g'(\rhon) \bigr) \, v .
  \label{minrho}
\Esist
In \eqref{minrho}, {\juerg we have $\tilde f=\tilde f_1+\tilde f_2$},
where $\tilde f_2$ is any smooth extension of $f_2$ to the whole of~$\erre$
and $\tilde f_1$ 
is the unique convex and lower semicontinuous extension of $f_1$
that satisfies $\tilde f(r)=+\infty$ if $r\not\in[0,1]$.
By the way, it is understood that the corresponding integral that appears in~\eqref{minrho}
is infinite if $\tilde f(v)$ does not belong to~$\Luno$.
Therefore, both functionals are well-defined and proper whenever $\mun\in V$ and $\rhon\in W$
({\paolo and} this implies boundedness of $g(\rhon)$ and~$g'(\rhon)$).
Moreover, {\paolo in view of} the above remarks, $J_1$ is continuous, and $J_2$ lower semicontinuous, on~$V$.

Now, we observe that equations \eqref{dprimabis} and \eqref{dsecondabis}, {\paolo when}
complemented by the regularity requirements in~\eqref{regdsoluz}
(which yield the homogeneous Neumann boundary conditions),
are the strong forms of the Euler-Lagrange variational equations 
{\paolo for} the stationary points of $J_1$ and~$J_2$, respectively.
More precisely, the strong form \eqref{dprimabis}
follows from the variational formulation thanks to the regularity theory of elliptic equations.
As far as \eqref{dsecondabis} is concerned,
the function $f'$ should be replaced {\paolo -- in principle, at least --} by the sum $\partial\tilde f_1+\tilde f_2'$,
where $\partial\tilde f_1$ is the subdifferential of~$\tilde f_1$.
However, once an $\Ldue$-estimate is obtained for the subdifferential
(and~standard arguments of the theory of maximal monotone operators 
(see, e.g., \cite{Brezis}), easily yield such an estimate),
the variational Euler-Lagrange equation can be written exactly {\juerg in the form} \eqref{dsecondabis},
{\paolo because} $\partial\tilde f_1$ is single-valued due to our assumptions on~$f_1$
(see \eqref{hpfp}, in particular).
{\juerg Consequently}, existence and uniqueness of the solution $(\mun,\rhon)$
follow if the functionals \eqref{minmu} and \eqref{minrho} are convex{\paolo , so that}
each of the correponding minimum problems has a unique solution.
This is {\paolo granted} for the first problem: 
indeed, $J_1$~is strictly convex and coercive{\paolo , because} $g$ is nonnegative.
The same {\paolo holds} for $J_2$, provided that the second derivative
of function $r\mapsto r^2/2+\tau f_2(r)$ is strictly positive on~$[0,1]$,
{\paolo which is the case} if $\tau\sup|f_2''|<1$.

It remains to prove that $\munp\geq0$.
To this {\juerg end}, we multiply \eqref{dprimabis} by $-\munp^-$,
where  $v^-=\max\{-v,0\}$ {\paolo denotes} the negative part of~$v$,
and integrate over~$\Omega$.
We obtain:
\Beq
  \iO \bigl( 1+g(\rhon)+g(\rhonp) \bigr) |\munp^-|^2
  + \tau \iO |\nabla\munp^-|^2
  = - \iO \bigl( 1+2g(\rhon) \bigr) \mun \munp^-
  \leq 0,
  \non
\Eeq
{\paolo because both} $g$ and $\mun$ are nonnegative.
This implies that $\munp^-=0$, {\paolo and hence that} $\munp\geq0$.


\section{Convergence}
\label{Convergenceproof}
\setcounter{equation}{0}

In this section, we prove Theorem~\ref{Convergenza}.
For convenience, we introduce one more vector, $(\xin)_{n=0}^N$,
and recall the definition of~$(\gn)_{n=0}^N$:
\Beq
  \xin := f_1'(\rhon)
  \aand
  \gn := g(\rhon)
  \quad \hbox{for $n=0,\dots,N$}.
  \label{defxingn}
\Eeq
Later on, we also use the interpolants of {\juerg these} vectors
according to Notation~\ref{Gennot}.
Our argument uses compactness and monotonicity methods.

\step
First a priori estimate

We multiply \eqref{dprimabis} by $\munp$ and integrate over~$\Omega$.
By accounting for the elementary identity~\eqref{elementare}, we obtain
\Bsist
  && \frac 12 \iO \munp^2 
  - \frac 12 \iO \mun^2
  + \frac 12 \iO |\munp-\mun|^2
  + \tau \iO |\nabla\munp|^2
  \non
  \\
  && + \iO \bigl(
    \gn \munp^2 + \gnp \munp^2 - 2 \gn \mun \munp
  \bigr)
  = 0 .
  \non
\Esist
As $\gn\munp^2+\gnp\munp^2-2\gn\mun\munp=\gnp\munp^2-\gn\mun^2+\gn(\munp-\mun)^2$,
we derive that
\Bsist
  && \iO \Bigl( \frac 12 + \gnp \Bigr) \munp^2 
  - \iO \Bigl( \frac 12 + \gn \Bigr) \mun^2
  + \iO \Bigl( \frac 12 + \gn \Bigr) |\munp-\mun|^2
  \non
  \\
  && \quad {}
  + \tau \iO |\nabla\munp|^2
  = 0 .
  \non
\Esist
{\paolo On} summing over $n=0,\dots,m-1$ with $1\leq m\leq N$, we conclude that
\Bsist
  && \iO \Bigl( \frac 12 + \gm \Bigr) \mum^2
  + \tau^2 \somma n0{m-1} \iO \Bigl( \frac 12 + \gn \Bigr) |\dtau\mun|^2
  + \tau \somma n0{m-1} \iO |\nabla\munp|^2
  \non
  \\
  && = \iO \Bigl( \frac 12 + \gamma_0 \Bigr) \muz^2 
  \quad \hbox{for $m=1,\dots,N$}.
  \non
\Esist
As  $g$ is nonnegative {\paolo and }hence $\gamma_i\geq0$, 
this implies that $\,\normaH\mum\leq c\,$
for $m=1,\dots,N$.
Thus, the above estimate also yields
\Beq
  \max_{m=1,\dots,N} \normaH\mum^2
  + \tau^2 \somma n0{N-1} \normaH{\dtau\mun}^2
  + \tau \somma n1N \normaV\mun^2
  \leq c .
  \label{primastimad}
\Eeq
In terms of the interpolants introduced in Notation~\ref{Gennot},
with the help of $\muz\in V$, \eqref{ouLinftyZ}--\eqref{ouLdueZ},
and \eqref{interpLinfty}--\eqref{interpLdue} we have that
\Bsist
  && \norma\overmu_{\L\infty H\cap\L2V}^2
  + \norma\undermu_{\L\infty H\cap\L2V}^2
  \non
  \\
  && \quad {}
  + \norma\hmu_{\L\infty H\cap\L2V}^2
  + \tau \norma{\dt\hmu}_{\L2H}^2
  \leq c.
  \label{primastima}
\Esist

\step
Second a priori estimate

In \eqref{dsecondabis}, we move $\rhon$ to the \lhs.
Then, we multiply by $\rhonp-\rhon$ and integrate over~$\Omega$.
We obtain
\Bsist
  && \iO |\rhonp-\rhon|^2
  + \frac \tau 2 \iO |\nabla\rhonp|^2
  - \frac \tau 2 \iO |\nabla\rhon|^2
  + \frac \tau 2 \iO |\nabla\rhonp - \nabla\rhon|^2
  \non
  \\
  && \quad {}
  + \tau \iO f'(\rhonp) (\rhonp-\rhon)
  = \tau \iO \mun g'(\rhon) (\rhonp-\rhon) .
  \label{perseconda}
\Esist
Now, we consider the last integral on the \lhs\ of~\eqref{perseconda}.
We split $f'=f_1'+f_2'$ and use the convexity assumption of~$f_1$
and boundedness for~$f_2'$.
We~get
\Beq
  \iO f'(\rhonp) (\rhonp-\rhon)
  \geq \iO \bigl( f_1(\rhonp) - f_1(\rhon) \bigr)
  - c \iO |\rhonp-\rhon|.
  \non
\Eeq
{\juerg Since also $g'$ is bounded, we infer from \eqref{perseconda} that}
\Bsist
  && \iO |\rhonp-\rhon|^2
  + \frac \tau 2 \iO |\nabla\rhonp|^2
  - \frac \tau 2 \iO |\nabla\rhon|^2
  + \frac \tau 2 \iO |\nabla\rhonp - \nabla\rhon|^2
  \non
  \\
  && \quad {}
  + \tau \iO \bigl( f_1(\rhonp) - f_1(\rhon) \bigr)
  \leq c \, \tau \iO (1+\mun) |\rhonp-\rhon| 
  \non
  \\
  && \leq \frac 12 \iO  |\rhonp-\rhon|^2 
  + c \, \tau^2 \iO (1+\mun^2)
  \leq \frac 12 \iO  |\rhonp-\rhon|^2 + c \, \tau^2,
  \non
\Esist
the last inequality being due to~\eqref{primastima}.
By dividing by~$\tau$, summing over $n=0,\dots,m-1$, 
and owing to the obvious inequality $m\tau\leq c$, we~conclude that
\Beq
  \tau \somma n0{m-1} \iO |\dtau\rhon|^2
  + \frac 12 \iO |\nabla\rhomp|^2 
  + \tau^2 \somma n0{m-1} \iO |\dtau(\nabla\rhon)|^2
  + \iO f_1(\rhomp)
  \leq c 
  \label{persecondabis}
\Eeq
for $m=0,\dots,N-1$.
As the term involving the difference quotient $\dtau\rhoNm$ is missing
in the first sum since $m\leq N-1$, we estimate it directly.
We multiply \eqref{dseconda}, written for $n=N-1$,
by $\tau\dtau\rhoNm$ and integrate over~$\Omega$.
We have
\Beq
  \tau \iO |\dtau\rhoNm|^2
  + \iO (\nabla\rhoN - \nabla\rhoNm) \cdot \nabla\rhoN
  + \iO f_1'(\rhoN) (\rhoN-\rhoNm)
  = \tau \iO \phi \, \dtau\rhoNm ,
  \non
\Eeq
where we have set
$\phi:=\mu_{N-1}\,g'(\rhoNm) - f_2'(\rhoNm)$.
{\ultime Owing to} the elementary identity \eqref{elementare}
and to the convexity of $f_1$ as before,
we have
\Bsist
  && \tau \normaH{\dtau\rhoNm}^2
  + \frac 12 \, \normaH{\nabla\rhoN}^2
  + \frac 12 \, \normaH{\nabla\rhoN-\nabla\rhoNm}^2
  + \iO f_1(\rhoN) 
  \non
  \\
  && \leq \frac 12 \, \normaH{\nabla\rhoNm}^2
  + \iO f_1(\rhoNm)
  + \tau \normaH\phi \, \normaH{\dtau\rhoNm}
  \non
  \\
  && \leq \frac 12 \, \normaH{\nabla\rhoNm}^2
  + \iO f_1(\rhoNm)
  + \frac \tau 2 \normaH\phi^2 
  + \frac \tau 2 \normaH{\dtau\rhoNm}.
  \non
\Esist
Now, we observe that the first two terms of the last line are bounded
by \eqref{persecondabis} written with $m=N-2$
and that $\phi$ is estimated in~$H$ thanks to~\eqref{primastimad}
and our assumptions of $g$ and~$f_2$.
Moreover, the last term of the first line can be ignored since $f_1$ is nonnegative.
Hence, we get the desired bound for the first term.
At this point, we can easily derive an estimate for $\normaH\rhom$ for $m=1,\dots,N$.
By using the obvious identity
$\rhom=\rhoz+\tau\somma n0{m-1}\drhon$
and the euclidean Schwarz and Young inequalities, we see~that
\Beq
  \normaH\rhom
  \leq \normaH\rhoz + \tau \somma n0{m-1} \normaH\drhon 
  \leq c + \frac \tau 2 \Bigl( m + \somma n0{m-1} \normaH\drhon^2 \Bigr)
  \leq c .
  \non
\Eeq
Hence, by recalling \eqref{persecondabis} and our last estimates,
we conclude that
\Beq
  \tau \somma n0{N-1} \normaH{\dtau\rhon}^2
  + \max_{m=1,\dots,N}\normaV\rhom^2
  + \tau^2 \somma n0{N-2} \normaH{\dtau(\nabla\rhon)}^2
  \leq c.
  \label{secondastimad}
\Eeq
In terms of the interpolants, \eqref{secondastimad} reads
(thanks also to $\rhoz\in V$ and to \eqref{ouLinftyZ}--\eqref{ouLdueZ} 
and \eqref{interpLinfty})
\Bsist
  \norma{\dt\hrho}_{\L2H}^2
  + \norma\overrho_{\L\infty V}^2
  + \norma\underrho_{\L\infty V}^2
  \qquad
  \non
  \\ 
  + \norma\hrho_{\L\infty V}^2
  + \tau \norma{\dt\nabla\hrho}_{L^2(0,T-\tau;H)}^2
  \leq c.
  \label{secondastima}
\Esist

\step
Third a priori estimate

We come back to~\eqref{dseconda} and rewrite it as (recall~\eqref{defxingn})
\Beq
  - \Delta\rhonp + \xinp
  = - \drhon + \mun g'(\rhon) - f_2'(\rhonp) .
  \non
\Eeq
Hence, a standard argument (multiplying by $-\Delta\rhonp$ and by~$\xinp$)
shows that the following estimate holds true
\Beq
  \normaH{\Delta\rhonp}
  + \normaH\xinp
  \leq c \normaH{- \drhon + \mun g'(\rhon) - f_2'(\rhonp)} \,.
  \non
\Eeq
Thus, we infer that
\Beq
  \normaH{\Delta\rhonp}^2
  + \normaH\xinp^2
  \leq c \bigl( \normaH\drhon^2 + \normaH\mun^2 + 1 \bigr) 
  \quad \hbox{for $n=0,\dots,N-1$}.
  \label{stimaDeltan}
\Eeq
Moreover, by using the regularity theory of elliptic equations,
we deduce that
\Beq
  \normaW\rhonp^2
  + \normaH\xinp^2
  \leq c \bigl( \normaV\rhonp^2 + \normaH\drhon^2 + \normaH\mun^2 + 1 \bigr) .
  \label{stimaWn}
\Eeq
Now, we multiply \eqref{stimaDeltan} by $\tau$ and sum over $n=0,\dots,m-1$.
By accounting for \eqref{primastimad} and~\eqref{secondastimad},
we conclude that 
\Beq
  \tau \somma n0{N-1} \normaW\rhonp^2
  + \tau \somma n0{N-1} \normaH\xinp^2
  \leq c .
  \label{terzastimad}
\Eeq
In terms of the interpolants, \eqref{terzastimad} yields
(by accounting for $\rhoz\in W$)
\Beq
  \norma\overrho_{\L2W}^2
  + \norma\underrho_{\L2W}^2
  + \norma\hrho_{\L2W}^2
  \leq c
  \label{terzastima}
\Eeq
besides an estimate for, e.g., $\norma\overxi$ in~${\L2H}$.

\step
Fourth a priori estimate

We write \eqref{dseconda} with $n+1$ in place of $n$
and take the difference between the obtained equality and \eqref{dseconda} itself.
Then we multiply {\juerg this} difference by $\dtau\rhonp$ and integrate over~$\Omega$.
We have
\Bsist
  && \iO (\dtau\rhonp - \dtau\rhon) \dtau\rhonp
  + \iO (\nabla\rhonpp - \nabla\rhonp) \cdot \dtau\nabla\rhonp
  \non
  \\
  && \quad {}
  + \iO \bigl( f'(\rhonpp) - f'(\rhonp) \bigr) \dtau\rhonp
  = \iO \bigl( \munp g'(\rhonp) - \mun g'(\rhon) \bigr) \dtau\rhonp .
  \qquad\qquad
  \label{perquarta}
\Esist
By accounting for the elementary identity~\eqref{elementare}, we get
\Beq
  \iO (\dtau\rhonp - \dtau\rhon) \dtau\rhonp
  = \frac 12 \iO |\dtau\rhonp|^2
  - \frac 12 \iO |\dtau\rhon|^2
  + \frac 12 \iO |\dtau\rhonp-\dtau\rhon|^2 .
  \non
\Eeq
Moreover, the second integral on the \lhs\ of \eqref{perquarta}
can be written in terms of $\dtau\rhonp$ in an obvious way.
Finally, by splitting $f'$ into $f_1'+f_2'$,
observing that the contribution due to the terms involving $f_1'$
is nonnegative since $f_1'$ is monotone
and moving the other ones to the \rhs,
we see that \eqref{perquarta} yields the inequality
\Bsist
  && \frac 12 \iO |\dtau\rhonp|^2
  - \frac 12 \iO |\dtau\rhon|^2
  + \frac 12 \iO |\dtau\rhonp-\dtau\rhon|^2
  + \tau \iO |\dtau\nabla\rhonp|^2
  \non
  \\
  && \leq - \iO \bigl( f_2'(\rhonpp) - f_2'(\rhonp) \bigr) \dtau\rhonp 
  + \iO \bigl( \munp g'(\rhonp) - \mun g'(\rhon) \bigr) \dtau\rhonp .
  \qquad
  \label{perquartabis}
\Esist
The first term on the \rhs\ of \eqref{perquartabis}
is easily treated {\juerg in the following way:}
\Beq
  - \iO \bigl( f_2'(\rhonpp) - f_2'(\rhonp) \bigr) \dtau\rhonp 
  \leq c \, \tau \iO |\dtau\rhonp|^2 .
  \label{perquartaF}
\Eeq
On the other hand, we have
\Bsist
  && \iO \bigl( \munp g'(\rhonp) - \mun g'(\rhon) \bigr) \dtau\rhonp
  \non
  \\
  && = \iO \munp \bigl( g'(\rhonp) - g'(\rhon) \bigr) \dtau\rhonp
  + \iO \bigl( \munp - \mun \bigr) g'(\rhon) \dtau\rhonp
  \non
  \\
  && \leq c \, \tau \iO \munp |\dtau\rhon| \, |\dtau\rhonp|
  + \iO \bigl( \munp - \mun \bigr) g'(\rhon) \dtau\rhonp \,.
  \non
\Esist
{\juerg Next,} we deal with the last integral by using equation~\eqref{dprima}.
{\ultime Owing to} our assumptions on~$g$, we obtain
\Bsist
  && \iO \bigl( \munp - \mun \bigr) g'(\rhon) \dtau\rhonp
  = - \tau \iO \frac {g'(\rhon)} {1+2\gn} \, 
    \bigl( \munp \, \dtau\gn - \Delta\munp
    \bigr) \dtau\rhonp
  \non
  \\
  && \leq c \, \tau \iO \munp \, |\dtau\rhon| \, |\dtau\rhonp|
  - \tau \iO \nabla\munp \cdot \nabla \Bigl( \frac {g'(\rhon)} {1+2\gn} \, \dtau\rhonp \Bigr)
  \non
  \\
  && \leq c \, \tau \iO \munp \, |\dtau\rhon| \, |\dtau\rhonp|
  \non
  \\
  && \quad {}
  - \tau \iO \frac {g'(\rhon)} {1+2\gn} \, \nabla\munp \cdot \nabla\dtau\rhonp 
  - \tau \iO \dtau\rhonp \, \nabla\munp \cdot \nabla \frac {g'(\rhon)} {1+2\gn}\,. 
  \non
\Esist
{\juerg We treat the last three terms separately}.
Thanks to the \holder, Sobolev, and Young inequalities,
and our assumption on~$g$, we have for every $\sigma>0$
\Bsist
  && \tau \iO \munp \, |\dtau\rhon| \, |\dtau\rhonp|
  \leq c \, {\gianni\tau} \, \norma\munp_4 \, \norma{\dtau\rhon}_{{\ultime H}} \, \norma{\dtau\rhonp}_4 
  \non
  \\
  && \leq \sigma \tau \normaV{\dtau\rhonp}^2
  + \frac {c \, \tau} \sigma \, \normaV\munp^2 \, \norma{\dtau\rhon}_{{\ultime H}}^2 
  \label{perquartaA}
  \\
  && - \tau \iO \frac {g'(\rhon)} {1+2\gn} \, \nabla\munp \cdot \nabla\dtau\rhonp
  \leq c \, \tau \iO |\nabla\munp| \, |\nabla\dtau\rhonp|
  \non
  \\ 
  && \leq \sigma \, \tau \normaH{\nabla\dtau\rhonp}^2
  + \frac {c \, \tau} \sigma \normaH{\nabla\munp}^2 
  \label{perquartaB}
  \\
  && - \tau \iO \dtau\rhonp \, \nabla\munp \cdot \nabla \frac {g'(\rhon)} {1+2\gn} 
  \leq c \, \tau \iO |\dtau\rhonp| \, |\nabla\munp| \, |\nabla\rhon|
  \non
  \\
  && \leq c \, \tau \norma{\dtau\rhonp}_4 \, \norma{\nabla\munp}_{{\ultime H}} \, \norma{\nabla\rhon}_4
  \vphantom\int
  \non
  \\
  && \leq \sigma \tau \normaV{\dtau\rhonp}^2
  + \frac {c\,\tau} \sigma \, \norma{\nabla\munp}_{{\ultime H}}^2 \, \bigl( \normaH\rhon^2 + \normaH{\Delta\rhon}^2 \bigr) \,.
  \label{perquartaC}
\Esist
Now, we rewrite \eqref{stimaDeltan} as
\Beq
  \normaH{\Delta\rhon}^2
  + {\gianni\normaH\xin^2}
  \leq c \bigl( \normaH{\dtau\rhonm}^2 + \normaH{\mu_{n-1}}^2 + 1 \bigr) 
  \quad \hbox{for $n=1,\dots,N$},
  \non
\Eeq
and note that we can allow the choice $n=0$
provided that we define 
\Beq
  \rho_{-1}:=\rhoz \ \hbox{ and, e.g., } \  \mu_{-1}:=0.
  \non
\Eeq
Hence, we can improve~\eqref{perquartaC}.
By using \eqref{primastima} and \eqref{secondastima} as well, we have
\Bsist
  && - \tau \iO \dtau\rhonp \, \nabla\munp \cdot \nabla \frac {g'(\rhon)} {1+2\gn}
  \non
  \\
  && \leq \sigma \tau \normaV{\dtau\rhonp}^2
  + \frac {c\,\tau} \sigma \, \norma{\nabla\munp}_{{\ultime H}}^2 \, \bigl( \normaH\rhon^2 + \normaH{\dtau\rhonm}^2 + \normaH{\mu_{n-1}}^2 + 1 \bigr) 
  \non
  \\
  && \leq \sigma \tau \normaV{\dtau\rhonp}^2
  + \frac {c\,\tau} \sigma \, \norma{\nabla\munp}_{{\ultime H}}^2 \, \bigl( \normaH{\dtau\rhonm}^2 + 1 \bigr) \,.
  \label{perquartaCbis}
\Esist
By recalling all these estimates, we see that \eqref{perquartabis} yields
\Bsist
  && \frac 12 \iO |\dtau\rhonp|^2
  - \frac 12 \iO |\dtau\rhon|^2
  + \frac 12 \iO |\dtau\rhonp-\dtau\rhon|^2
  + \tau \iO |\dtau\nabla\rhonp|^2
  \non
  \\
  && \leq c \, \tau \iO |\dtau\rhonp|^2
  + 3 \sigma \tau \normaV{\dtau\rhonp}^2
  + \frac {c \, \tau} \sigma \, \normaV\munp^2 \, \normaH{\dtau\rhon}^2 
  \non
  \non
  \\
  && \quad {}
  + \frac {c\,\tau} \sigma \, \normaH{\nabla\munp}^2 \, \bigl( \normaH{\dtau\rhonm}^2 + 1 \bigr) \,.
  \non
\Esist
Now, just by changing the value of the constant $c$ in front of the first integral on the \rhs,
we can replace the last integral on the \lhs\ by $\normaV{\dtau\rhonp}^2$.
Then, we choose $\sigma=1/4$ and rearrange. We obtain
\Bsist
  && \frac 12 \iO |\dtau\rhonp|^2
  - \frac 12 \iO |\dtau\rhon|^2
  + \frac 12 \iO |\dtau\rhonp-\dtau\rhon|^2
  + \frac \tau 4 \, \normaV{\dtau\rhonp}^2
  \non
  \\
  && \leq c \, \tau \iO |\dtau\rhonp|^2
  + c \, \tau \, \normaV\munp^2 \, \normaH{\dtau\rhon}^2 
  + c \,\tau \, \normaH{\nabla\munp}^2 \, \bigl( \normaH{\dtau\rhonm}^2 + 1 \bigr) \,.
  \non
\Esist
At this point, by assuming $m\leq N-1$, we sum over $n=0,\dots,m-1$ and have
\Bsist
  && \frac 12 \iO |\dtau\rhom|^2
  + \frac 12 \, \somma n0{m-1} \iO |\dtau\rhonp-\dtau\rhon|^2
  + \frac \tau 4 \, \somma n0{m-1} \normaV{\dtau\rhonp}^2
  \non
  \\
  && \leq \frac 12 \iO |\dtau\rhoz|^2
  + c \, \tau \somma n0{m-1} \iO |\dtau\rhonp|^2
  + c \, \tau \somma n0{m-1} \normaV\munp^2 \, \normaH{\dtau\rhon}^2 
  \non
  \\
  && \quad {}
  + c\,\tau \somma n0{m-1} \normaH{\nabla\munp}^2 \, \normaH{\dtau\rhonm}^2 
  + c\,\tau \somma n0{m-1} \normaH{\nabla\munp}^2 \,.
  \label{perquartater}
\Esist
The second and {\paolo the} last terms on the \rhs\ of \eqref{perquartater}
have been already estimated 
by \eqref{secondastimad} and~\eqref{primastimad}, respectively.
To treat the first term, we write \eqref{dseconda} with $n=0$ 
and add $\Delta\rhoz$ to both sides.
Then, we multiply the resulting equality by $\dtau\rhoz$ and integrate over~$\Omega$.
{\paolo After a rearrangement,} owing to \eqref{hpfd} and the assumptions on the initial data
(see~\eqref{hprhoz}, in particular), we obtain:
\Bsist
  &&\iO |\dtau\rhoz|^2
  + \tau \iO |\nabla\dtau\rhoz|^2
  + \iO \bigl( f'_1 (\rho_1) - f'_1 (\rhoz) \bigr) \dtau\rhoz 
  \non 
  \\
   &&{}= \iO \bigl( \Delta\rhoz + \muz g'(\rhoz) - f'_1 (\rhoz) - f_2'(\rho_1) \bigr) \dtau\rhoz 
  \leq c \normaH{\dtau\rhoz} . 
  \non
\Esist
As  $ \bigl( f'_1 (\rho_1) - f'_1 (\rhoz) \bigr) \dtau\rhoz \geq 0 $ due to the monotonicity of  $f'_1,$ we immediately deduce~that
\Beq
  \normaH{\dtau\rhoz}^2 + \tau \normaH{\nabla\dtau\rhoz}^2 \leq c .
  \label{stimadtrhoz}
\Eeq
In particular, the desired estimate for~$\normaH{\dtau\rhoz}$ is achieved.
Therefore, {\paolo on} recalling that $\dtau\rho_{-1}=0$ {\paolo because} $\rho_{-1}=\rhoz$,
we see that \eqref{perquartater} yields:
\Bsist
  && \normaH{\dtau\rhom}^2
  + \somma n0{m-1} \normaH{\dtau\rhonp-\dtau\rhon}^2
  + \tau \somma n0{m-1} \normaV{\dtau\rhonp}^2
  \non
  \\
  \noalign{\allowbreak}
  && \leq c
  + c \, \tau \somma n0{m-1} \normaV\munp^2 \, \normaH{\dtau\rhon}^2 
  + c\,\tau \somma n1{m-1} \normaH{\nabla\munp}^2 \, \normaH{\dtau\rhonm}^2 
  \non
  \\
  && \leq c
  + c \, \tau \somma n0{m-1} \normaV\munp^2 \, \normaH{\dtau\rhon}^2 
  + c\,\tau \somma n0{m-2} \normaH{\nabla\mu_{n+2}}^2 \, \normaH{\dtau\rhon}^2 
  \non 
  \\
  && \leq C_1
  + C_2 \, \tau \somma n0{m-1} \bigl( \normaV\munp^2 + \normaV{\mu_{n+2}}^2 \bigr) \, \normaH{\dtau\rhon}^2
  \non 
\Esist
for $m=0,\dots,N-1$. 
Hence, we can apply the discrete Gronwall lemma 
(see~\eqref{dgronwall}, where $N$ is {\paolo to be replaced here} by $N-1$)
and deduce~that
\Bsist
  && \normaH{\dtau\rhom}^2
  + \somma n0{m-1} \normaH{\dtau\rhonp-\dtau\rhon}^2
  + \tau \somma n0{m-1} \normaV{\dtau\rhonp}^2
  \non
  \\
  && \leq C_1 \exp \Bigl( C_2 \tau \somma n0{m-1} \bigl( \normaV\munp^2 + \normaV{\mu_{n+2}}^2 \bigr) \Bigr)
  \non
\Esist
for $m=0,\dots,N-1$.
{\paolo O}wing to \eqref{primastimad}, we infer that
\Beq
  \normaH{\dtau\rhom}^2
  + \somma n0{m-1} \normaH{\dtau\rhonp-\dtau\rhon}^2
  + \tau \somma n0{m-1} \normaV{\dtau\rhonp}^2
  \leq c
  \quad \hbox{for $m=0,\dots,N-1$};
  \non
\Eeq
{\paolo moreover, on} using the estimates of $\dtau\rhoz$ and $\nabla\dtau\rhoz$ given by~\eqref{stimadtrhoz},
we conclude~that
\Beq
  \max_{m=0,\dots,N-1} \normaH{\dtau\rhom}^2
  + \somma n0{N-1} \normaH{\dtau\rhon-\dtau\rhonm}^2
  + \tau \somma n0{N-1} \normaV{\dtau\rhon}^2
  \leq c .
  \label{quartastimad}
\Eeq
In particular, \eqref{quartastimad} yields:
\Beq
  \norma{\dt\hrho}_{\L\infty H\cap\L2V} \leq c.
  \label{quartastima}
\Eeq

\step
Fifth a priori estimate

We improve \accorpa{terzastimad}{terzastima}.
{\paolo O}wing to \eqref{stimaWn}, on using~\eqref{quartastimad}
in addition to previous estimates,
we immediately obtain (cf.~also \eqref{hprhoz}) {\juerg that}
\Bsist
  && \normaW\rhom + \normaH\xim 
  \leq c 
  \quad \hbox{for $m=0,\dots,N$,}
  \label{quintastimad}
  \\
  && \norma\overrho_{\L\infty W}^2
  + \norma\underrho_{\L\infty W}^2
  + \norma\hrho_{\L\infty W}^2
  \leq c \,,
  \vphantom\int
  \label{quintastima}
\Esist
{\juerg as well as} an estimate for, e.g., $\overxi$ in $\L\infty H$.

\step
Sixth a priori estimate

We rewrite \eqref{dprima} in the form
\Beq
  \bigl( 1+\gn+\gnp \bigr) \dmun
  - \Delta\munp 
  = - \mun \, \dgn \,.
  \non
\Eeq
{\paolo W}e test {\juerg this} equality by $(\munp-\mun)$, and integrate over~$\Omega$.
We obtain
\Beq
  \tau \iO \bigl( 1+\gn+\gnp \bigr) |\dmun|^2
  + \iO (\nabla\munp-\nabla\mun) \cdot \nabla\munp
  = - \tau \iO \mun \, \dgn \, \dmun .
  \non
\Eeq
As $g$ is nonnegative and Lipschitz continuous, we infer that
\Bsist
  && \tau \iO |\dmun|^2
  + \frac 12 \iO |\nabla\munp|^2 
  - \frac 12 \iO |\nabla\mun|^2 
  + \frac 12 \iO |\nabla(\munp-\mun)|^2 
  \non
  \\
  && \leq c \, \tau \iO \mun \, |\drhon| \, |\dmun|
  \leq c \, \tau \, \norma\mun_4 \, \norma\drhon_4 \, \norma\dmun_2
  \non
  \\
  && \leq \frac \tau 2 \, \normaH\dmun^2
  + c \, \tau \normaV\drhon^2 \, \bigl( \normaH{\nabla\mun}^2 + \normaH\mun^2 \bigr)
  \non
  \\
  && \leq \frac \tau 2 \, \normaH\dmun^2
  + c \, \tau \normaV\drhon^2 \, \normaH{\nabla\mun}^2 + c \, \tau \normaV\drhon^2\,,
  \non
\Esist
the last inequality being due to~\eqref{primastimad}.
By rearranging and summing over $n=0,\dots,m-1$ with $1\leq m\leq N$, we~get:
\Bsist
  && \frac \tau 2 \, \somma n0{m-1} \normaH\dmun^2
  + \frac 12 \, \normaH{\nabla\mum}^2
  + \frac {\tau^2} 2 \, \somma n0{m-1} \normaH{\nabla\dmun}^2
  \non
  \\  
  && \leq \frac 12 \, \normaH{\nabla\muz}^2
  + c \, \tau \somma n0{m-1} \normaV\drhon^2 \, \normaH{\nabla\mun}^2 
  + c \, \tau \somma n0{m-1} \normaV\drhon^2 
  \non
  \\  
  && \leq c
  + c \, \tau \somma n0{m-1} \normaV\drhon^2 \, \normaH{\nabla\mun}^2 \,,
  \non
\Esist
where we have used~\eqref{quartastimad}.
Now, we first apply the discrete Gronwall lemma~\eqref{dgronwall}
and then account for~\eqref{quartastimad} once more.
We obtain, for $m=1,\dots,N$,
\Beq
  \tau \somma n0{m-1} \normaH\dmun^2
  + \normaH{\nabla\mum}^2
  + \tau^2 \, \somma n0{m-1} \normaH{\nabla\dmun}^2
  \leq c.
  \label{persesta}
\Eeq
Next, by \eqref{dprima}, the \holder\ and Sobolev inequalities 
and the Lipschitz continuity of~$g$, we infer~that
\Bsist
  && \normaH{\Delta\munp}
  \leq c \bigl( \normaH\dmun + \normaH{\munp\dgn} \bigr)
  \leq c \bigl( \normaH\dmun + \norma\munp_4 \, \norma\dgn_4 \bigr)
  \non
  \\
  && \leq c \bigl( \normaH\dmun + \norma\munp_4 \, \norma\drhon_4 \bigr)
  \leq c \bigl( \normaH\dmun + \normaV\munp \, \normaV\drhon \bigr);
  \non
\Esist
{\paolo note that in the last product} we can ignore the factor $\normaV\munp$,
due to \eqref{primastima} and~\eqref{persesta},
provided we update the last value of~$c$.
By squaring, summing up, and multiplying by~$\tau$, we thus obtain for $m=1,\dots,N$ {\juerg the estimate}
\Beq
  \tau \somma n0{m-1} \normaH{\Delta\munp}^2
  \leq c \, \tau \somma n0{m-1} \normaH\dmun^2 
  + c \, \tau \somma n0{m-1} \normaV\drhon^2\,,
  \non 
\Eeq
and we can replace the $H$-norm of $\Delta\munp$ 
by the $W$-norm of~$\munp$ thanks to~\eqref{primastimad}.
We collect this and \eqref{persesta} and account for \eqref{quartastimad} and $\muz\in V$.
We have:
\Beq
  \tau \somma n0{N-1} \normaH\dmun^2
  + \max_{n=0,\dots,N} \normaH{\nabla\mun}^2
  + \tau \somma n0{N-1} \normaW\munp^2
  \leq c \,,
  \label{sestastimad}
\Eeq
so that 
\Beq
  \norma{\dt\hmu}_{\L2H}
  + \norma\overmu_{\L\infty V\cap\L2W}
  + \norma\undermu_{\L\infty V}
  + \norma\hmu_{\L\infty V}
  \leq c .
  \label{sestastima}
\Eeq
We note that \eqref{persesta} also gives the non-sharp estimate
\Beq
  \tau \norma{\nabla\dt\hmu}_{\L2H}^2 \leq c .
  \label{nonsharp}
\Eeq

\step
Limit and conclusion

By standard weak compactness {\juerg results}, we find some convergent subsequence
for the interpolants.
{\paolo Therefore, in principle,} it is understood that the convergence {\juerg that} we refer to
holds for a subsequence.
However, once we prove that the limit we find is the solution $(\mu,\rho)$
to problem \Pbl, then the whole family of interpolants is convergent,
due to uniqueness.
For the reader's convenience, we select some estimates 
among those we have proved in the previous steps.
{\paolo These are:}
\Bsist
  && \norma\overmu_{\L\infty V\cap\L2W}
  + \norma\undermu_{\L\infty V}
  + \norma\hmu_{\L\infty V}
  \leq c
  \label{stimemu}
  \\
  && \norma\overrho_{\L\infty W}
  + \norma\underrho_{\L\infty W}
  + \norma\hrho_{\L\infty W}
  \leq c
  \label{stimerho}
  \\
  && \norma{\dt\hmu}_{\L2H}
  + \norma{\dt\hrho}_{\L\infty H\cap\L2V} \leq c.
  \label{stimedthat}
\Esist
Now, we observe that \eqref{stimedthat} and 
\eqref{interpLdue} imply that
\Bsist
  && \norma{\overmu-\hmu}_{\L2H}
  + \norma{\undermu-\hmu}_{\L2H}
  \leq c \, \tau
  \label{stimadiffmu}
  \\
  && \norma{\overrho-\hrho}_{\L2V}
  + \norma{\underrho-\hrho}_{\L2V}
  \leq c \, \tau .
  \label{stimadiffrho}
\Esist
This yields, in particular, that the weak limits we find for 
$\overmu$, $\undermu$, and $\hmu$, 
by using \eqref{stimemu} and weak compactness results coincide
and that the same happens {\paolo for} $\overrho$, $\underrho$, and~$\hrho$.
Therefore, we can conclude that some functions $\mu$ and $\rho$ exist
such~that
\Bsist
  & \overmu, \, \undermu, \, \hmu \to \mu
  & \quad \hbox{weakly star in $\L\infty V$,}
  \label{convmu} 
  \\
  & \overmu \to \mu
  & \quad \hbox{weakly in $\L2W$,}
  \label{convmubis} 
  \\
  & \overrho, \, \underrho, \, \hrho \to \rho
  & \quad \hbox{weakly star in $\L\infty W$,}
  \label{convrho} 
  \\
  & \dt\hmu \to \dt\mu
  & \quad \hbox{weakly in $\L2H$,}
  \label{convdtmu} 
  \\
  & \dt\hrho\to \dt\rho
  & \quad \hbox{weakly star in $\L\infty H\cap\L2V$}\,.
  \label{convdtrho} 
\Esist
{\juerg Now} we prove that $(\mu,\rho)$ satisfies \Regsoluz\ and solves problem \Pbl.

We remark that the topology {\paolo alluded {\pier to}} in the statement
of Theorem~\ref{Convergenza} is precisely the topology
associated {\paolo with} the convergence{\paolo s specified in} \accorpa{convmu}{convdtrho}.
Clearly, \accorpa{regmu}{mupos} are fulfilled.
Moreover, the Cauchy conditions~\eqref{cauchy} are satisfied,
{\paolo because} $(\hmu,\hrho)$ converges to $(\mu,\rho)$ at least weakly in~$\C0H$.
Therefore, it remains to check that \eqref{regfprho} holds 
and that equations \accorpa{prima}{seconda} are satisfied.
To do that, we read the discrete problem \Discreto\ 
in terms of the interpolants.
We have:
\Bsist
  && \bigl( 1+2\underg \bigr) \dt\hmu + \overmu \dt\hg - \Delta\overmu = 0,
  \label{discrmu}
  \\
  && \dt\hrho - \Delta\overrho + f'(\overrho) = \undermu \, g'(\underrho) .
  \label{discrrho}
\Esist
Hence, the main problem consists in identifying correctly the limits of the 
nonlinear terms and those of the products.
To this {\juerg end}, we recover some strong convergence
(without looking for sharpness, since it is not necessary).
We~first recall that the embeddings $V\subset H$ and $W\subset\Cx0$ are compact,
{\paolo so that} we can apply~\cite[Sect.~8, Cor.~4]{Simon} and deduce~that
\Bsist
  && \hmu \to \mu
  \quad \hbox{strongly in $\C0H$ and \aeQ},
  \label{strongmu}
  \\
  && \hrho \to \rho
  \quad \hbox{strongly in $\C0{\Cx0}=C^0(\overline Q)$}.
  \label{strongrho}
\Esist
By combining this with \eqref{stimadiffmu} and~\eqref{stimadiffrho},
we infer that
\Beq
  \overmu,\, \undermu \to \mu
  \aand
  \overrho, \, \underrho \to \rho
  \quad \hbox{strongly in $\L2H$ and \aeQ}.
  \qquad
  \label{strongbis}
\Eeq
We point out that {\paolo a.e.-}convergence actually holds 
for a subsequence. 
{\corr
As $f_2'$, $g$, and $g'$ are Lipschitz continuous on~$[0,1]$,
we deduce~that
\Beq
  \phi(\overrho), \, \phi(\underrho) \to \phi(\rho)
  \quad \hbox{strongly in $\L2H$\quad for $\phi=f_2',g,g'$}.
  \label{strongnonlin}
\Eeq
On the other hand, by comparison in~\eqref{discrrho},
we see that $f_1'(\overrho)$ remains bounded in $\L\infty H$,
so that $f_1'(\overrho)$ converges (for a subsequence) to~some $\xi$ 
in the weak star topology of such a space. As $f_1'$ induces a maximal monotone operator 
on $L^2(Q)\times L^2 (Q)$ (cf., e.g., \cite[Exemple~2.3.3, p.~25]{Brezis}),
$f_1'(\overrho) \to \xi$  and $ \overrho \to \rho $ weakly in $L^2(Q)$, 
and 
$$ \limsup_{\tau\searrow 0} \int_Q  f_1'(\overrho) \overrho \leq   \int_Q  \xi \rho, $$
owing to standard results in the theory of maximal monotone 
operators (one may see \cite[Prop.~2.5, p.~27]{Brezis}),
we deduce that $0<\rho<1$ and $\xi=f_1'(\rho)$ \aeQ.
In particular, \eqref{regfprho}~holds.
Furthermore, we also have
\Beq
  |\dt\hg|
  = |\dtau\gn|
  = |\dtau g(\rhon)|
  \leq c |\dtau\rhon|
  = c |\dt\hrho|
  \quad \hbox{a.e.\ in $I_{n+1}$, for $n=0,\dots,N-1$}\,,
  \non
\Eeq
so that \eqref{stimedthat} yields an estimate of $\dt\hg$ in~$\L\infty H$.
Hence, thanks to~\eqref{interpLinfty}, we~have
\Beq
  \norma{\hg - \overg}_{\L\infty H} 
  \leq c \, \tau \norma{\dt\hg}_{\L\infty H} 
  \leq c \, \tau\,,
  \non
\Eeq
whence even $\hg$ converges to $g(\rho)$, e.g., strongly in~$\LQ2$.
Then, we deduce~that
\Beq
  \dt\hg \to \dt g(\rho)
  \quad \hbox{weakly star in $\L\infty H$}.
  \label{convdthg}
\Eeq
Finally, as {\paolo to the limits of the products in \accorpa{discrmu}{discrrho}},
we can infer~that
\Beq
  \underg \dt\hmu \to g(\rho) \dt\mu , \quad
  \overmu \dt\hg \to \mu \dt g(\rho) , \quad
  \overmu g'(\underrho) \to \mu g'(\rho),
  \quad \hbox{weakly in $\LQ1$}.
  \non
\Eeq
Therefore, \accorpa{prima}{seconda} follow from \accorpa{discrmu}{discrrho},
and the proof is complete. In particular, let us stress that the so found pair $(\mu, \rho) $  
solves \eqref{regmu}--\eqref{cauchy} and then it must coincide with the unique solution 
$(\mu, \rho) $ of the continuous problem given by Theorem~\ref{Esistcont}.}

{\corr
As a by-product of the above proof, it turns out that
\Beq
  \rhominbis \leq \hrho,\, \overrho,\, \underrho \leq \rhomaxbis
  \quad \hbox{in $Q$,\quad for some $\rhominbis,\rhomaxbis\in(0,1),$}
  \label{farzeroone}
\Eeq
provided that $\tau$ is small enough. 
Indeed, take $\rhominbis\in(0,\rhomin)$ and $\rhomaxbis\in(\rhomax,1)$,
with $\rhomin,\rhomax\in(0,1)$ given by Theorem~\ref{Esistcont}.
That \eqref{farzeroone} holds for $\hrho$ follows 
from the uniform convergence given by~\eqref{strongrho}.
This means that the same bounds hold for $\rhon$, $n=0,\dots,N$
(where $(\rhon)_{n=0}^N$ is the vector associated with~$\hrho$),
and hence also for $\overrho$ and~$\underrho$.
}


\section{Proof of Theorem \ref{Errore}}
\label{ProofErrore}
\setcounter{equation}{0}

In this section, we prove Theorem~\ref{Errore}.
It is understood that $\tau$ is as small as {\paolo needed; oftentimes,
we do not pause and quantify such smallness precisely.}
First of all, we remind the reader that the interpolants $\hrho$, $\underrho$, and $\overrho$
are uniformly far for $0$ and~$1$
(see~\eqref{farzeroone} and the subsequent lines).
Therefore, {\paolo without loss of generality, we can assume that the derivative function} $f'$ is Lipschitz continuous.
We need {\juerg additional} a~priori estimates.

\step
Auxiliary a priori estimates

We prepare an estimate for $\normaH{\nabla\dtau\rhoz}$.
To this end, we notice that \eqref{dseconda} with $n=0$ can be written~as
\Beq
  \dtau\rhoz - \tau \Delta\dtau\rhoz
  =  f'(\rhoz) - f'(\rho_1) - \psi_0\,,
  \label{auxzero}
\Eeq
where $ \psi_0 := -\Delta\rhoz + f'(\rhoz) - \muz g'(\rhoz) .$
As $\psi_0\in V$ by~\eqref{dahperrore}, we can test \eqref{auxzero} by $-\Delta\dtau\rhoz$ and integrate by parts.  
{\paolo In view of} the Lipschitz continuity of~$f'$, we find out~that
\Bsist
  && \normaH{\nabla\dtau\rhoz}^2 
  + \tau \normaH{\Delta\dtau\rhoz}^2
  \leq c \, \tau \normaH{\dtau\rhoz} \, \normaH{\Delta\dtau\rhoz}
  + \normaH{\nabla\psi_0} \, \normaH{\nabla\dtau\rhoz}
  \non
  \\
  && \leq \frac \tau 2 \, \normaH{\Delta\dtau\rhoz}^2
  + c \, \tau \normaH{\dtau\rhoz}^2
  + \frac 12 \, \normaH{\nabla\dtau\rhoz}^2 
  + c .
  \non
\Esist
By accounting for \eqref{stimadtrhoz}, we obtain the desired estimate 
\Beq
  \normaH{\nabla\dtau\rhoz} \leq c .
  \label{nabladtrhoz}
\Eeq
Let us come now to the basic estimate we need.
We write \eqref{dseconda} with $(n+1)$ in place of $n$,
and take the difference between the {\paolo so-}obtained equality and \eqref{dseconda} itself.
Then, we multiply {\juerg this} difference by $-\Delta\drhonp$ and integrate over~$\Omega$.
We easily have, for $n=0,\dots,N-2$, {\paolo that}
\Bsist
  && \iO (\nabla\drhonp - \nabla\drhon) \cdot \nabla\drhonp
  + \iO (\Delta\rhonpp - \Delta\rhonp) \, \Delta\drhonp
  \non
  \\
  && 
  = - \iO \bigl( f'(\rhonpp) - f'(\rhonp) \bigr) (-\Delta\drhonp)
  \non
  \\
  && \quad {}
  + \iO \bigl( \munp g'(\rhonp) - \mun g'(\rhon) \bigr) (-\Delta\drhonp) .
  \label{peraux}
\Esist
By the elementary identity~\eqref{elementare}, 
the first integral is equal~to
\Beq
  \frac 12 \iO |\nabla\drhonp|^2
  - \frac 12 \iO |\nabla\drhon|^2
  + \frac 12 \iO |\nabla\drhonp-\nabla\drhon|^2 .
  \non
\Eeq
On the other hand, we obviously have {\paolo that}
\Beq
  \iO (\Delta\rhonpp - \Delta\rhonp) \, \Delta\drhonp
  = \tau \iO |\Delta\drhonp|^2 .
  \non
\Eeq
Now, we deal with the \rhs\ of~\eqref{peraux}.
By Lipschitz continuity, we deduce~that
\Bsist
  && - \iO \bigl( f'(\rhonpp) - f'(\rhonp) \bigr) (-\Delta\drhonp)
  \leq c \iO |\rhonpp-\rhonp| \, |\Delta\drhonp|
  \non
  \\
  && = c \, \tau \iO |\drhonp| \, |\Delta\drhonp|
  \leq \frac \tau 4 \iO |\Delta\drhonp|^2
  + c \, \tau \iO |\drhonp|^2 .
  \non
\Esist
As far as the last term of \eqref{peraux} is concerned,
we combine the above elementary argument with the \holder\ and Young inequalities 
and the Sobolev embedding $V\subset\Lx4$.
We~{\paolo find:}
\Bsist
  && \iO \bigl( \munp g'(\rhonp) - \mun g'(\rhon) \bigr) (-\Delta\drhonp)
  \non
  \\
  && \leq c \iO \bigl( \munp |\rhonp-\rhon| + |\munp-\mun| \bigr) |\Delta\drhonp|
  \non
  \\
  && \leq c \, \tau \bigl( \norma\munp_4 \norma\drhon_4 + \normaH\dmun \bigr) \normaH{\Delta\drhonp}
  \non
  \\
  && \leq \frac \tau 4 \iO |\Delta\drhonp|^2
  + c \tau \normaV\munp^2 \, \normaV\drhon^2  
  + c \tau \normaH\dmun^2 \, .
  \non
\Esist
By collecting the inequalities we have obtained,
we see that \eqref{peraux} yields:
\Bsist
  && \frac 12 \iO |\nabla\drhonp|^2
  - \frac 12 \iO |\nabla\drhon|^2
  + \frac 12 \iO |\nabla\drhonp-\nabla\drhon|^2
  + \frac \tau 2 \iO |\Delta\drhonp|^2  
  \non
  \\
  \non
  \\
  && \leq c \, \tau \iO |\drhonp|^2
  + c \tau \normaV\munp^2 \, \normaV\drhon^2  
  + c \tau \normaH\dmun^2 \, .
  \non
\Esist
At this point, we sum over $n=0,\dots,m-1$, with $1\leq m\leq N-1$, and deduce that
\Bsist
  && \frac 12 \iO |\nabla\dtau\rhom|^2
  + \frac 12 \, \somma n0{m-1} \iO |\nabla\drhonp-\nabla\drhon|^2
  + \frac \tau 2 \, \somma n0{m-1} \iO |\Delta\drhonp|^2  
  \non
  \\
  && \leq \frac 12 \iO |\nabla\dtau\rhoz|^2
  + c \, \tau \somma n0{N-2} \normaH\drhonp^2
  \non
  \\
  && \quad {}
  + c \max_{n=0,\dots,N-1} \normaV\munp^2 \,\, \tau \somma n0{N-1} \normaV\drhon^2 
  + c \, \tau \somma n0{N-1} \normaH\dmun^2 \, .
  \qquad
  \label{perauxbis}
\Esist
The first term on the \rhs\ of~\eqref{perauxbis} is estimated by~\eqref{nabladtrhoz}{\paolo ; all} other terms on the \rhs\ 
have been estimated already
(cf.\ \eqref{quartastima} and~\eqref{sestastima}).
Therefore, by recalling~also \eqref{quartastimad},
we conclude~that
\Bsist
  && \max_{m=0,\dots,N-1} \normaV{\drhom}^2
  + \tau \somma n0{N-2} \normaH{\Delta\drhonp}^2
  \leq c ,
  \label{stimaauxd}
  \\
  && \norma{\dt\hrho}_{\L\infty V}
  + \norma{\Delta\dt\hrho}_{\L2H}
  \leq c . 
  \vphantom\int
  \label{stimaaux}
\Esist

\step
Consequence

{\paolo In view of} the regularity theory for elliptic equations
and the continuous embedding $W\subset\Linfty$, we derive from~\eqref{stimaaux} {\paolo that}
\Beq
  \norma{\dt\hrho}_{\L2W} \leq c 
  \aand
  \norma{\dt\hrho}_{\L2\Linfty} \leq c .
  \label{daaux}
\Eeq
Moreover, as the second \eqref{daaux} means an estimate of the difference quotients
associated to the vector $(\rhon)_{n=0}^N$, and {\paolo as} $g$ is Lipschitz continuous,
a similar estimate holds for the vector~$(g(\rhon))_{n=0}^N$ {\paolo 
(see}~\eqref{defxingn}),
and we infer~that
\Beq
  \norma{\dt\hg}_{\L2\Linfty} \leq c .
  \label{daauxbis}
\Eeq
Furthermore, by applying \eqref{interpLdue}, we see that \eqref{stimaaux} also implies 
{\juerg that}
\Beq
  \norma{\Delta(\overrho-\hrho)}_{\L2H} \leq c \, \tau. 
  \label{daauxter}
\Eeq

\step
Proof of Theorem~\ref{Errore}

A possible strategy could be the following{\paolo : to} 
multiply the difference between \eqref{discrmu} and \eqref{prima}
by~$(\hmu-\mu)$, and the difference between \eqref{discrrho} and \eqref{seconda}
by~$\dt(\hrho-\rho)${\paolo ; then, to} sum up and start estimating.
However, in order to split {\paolo c}alculations
and give more transparence to the proof,
we prefer to proceed with {\paolo those} pairs of equation separately,
and collect the inequalities we obtain later on.
So, we first consider just one couple,
for instance, \eqref{discrrho} and~\eqref{seconda}.
We multiply their difference by~$\dt(\hrho-\rho)$,
integrate over~$Q_t$, where $t\in(0,T)$ is arbitrary,
and add the same integral to both {\juerg sides} for convenience.
We obtain:
\Bsist
  && \intQt |\dt(\hrho-\rho)|^2
  + \frac 12 \, \normaV{(\hrho-\rho)(t)}^2
  \non
  \\
  && = \intQt
  \Bigl\{
    - \Delta(\hrho-\overrho)
    - \bigl( f'(\overrho) - f'(\rho) \bigr)
  \non
  \\
  && \qquad {}
    + g'(\underrho) (\undermu-\mu)
    + \mu \bigl( g'(\underrho) - g'(\rho) \bigr)
    + (\hrho-\rho)
  \Bigr\} \dt(\hrho-\rho)
  \non
  \\
  && \leq \frac 12 \intQt |\dt(\hrho-\rho)|^2
  \non
  \\
  && \quad {}
  + c \intQt
  \Bigl\{
    |\Delta(\hrho-\overrho)|^2
    + |\overrho-\rho|^2
    + |\undermu-\mu|^2
    + |\underrho-\rho|^2
    + |\hrho-\rho|^2
  \Bigr\} .
  \qquad
  \label{pererrorerho}
\Esist
In the above inequality,
we have used the Lipschitz continuity of $f'$ and~$g'$, 
and the boundedness of~$\mu$.
Now, we estimate the last integral of~\eqref{pererrorerho}.
Thanks to~\eqref{daauxter}, we have~{\paolo that}
\Beq
  \intQt |\Delta(\hrho-\overrho)|^2 \leq c \, \tau^2.
  \non
\Eeq
On the other hand, owing to~\eqref{stimadiffrho}, we obtain:
\Bsist
  && \intQt \bigl( |\overrho-\rho|^2 + |\underrho-\rho|^2 + |\hrho-\rho|^2 \bigr)
  \leq c \intQt \bigl( |\overrho-\hrho|^2 + |\underrho-\hrho|^2 + |\hrho-\rho|^2 \bigr)
  \non
  \\
  && \leq c \, \tau^2 + c \intQt |\hrho-\rho|^2 .
  \non
\Esist
Similarly, we have, by~\eqref{stimadiffmu}, that
\Beq
  \intQt |\undermu-\mu|^2
  \leq c \intQt \bigl( |\undermu-\hmu|^2 + |\hmu-\mu|^2 \bigr)
  \leq c \, \tau^2 
  + c \intQt |\hmu-\mu|^2 .
  \non
\Eeq
By collecting the above inequalities,
we see that \eqref{pererrorerho} and the Gronwall lemma yield:
\Beq
  \intQt |\dt(\hrho-\rho)|^2
  + \normaV{(\hrho-\rho)(t)}^2
  \leq c
  \Bigl\{ 
    \tau^2  + \intQt |\hrho-\rho|^2
    + \intQt |\hmu-\mu|^2 
  \Bigr\} 
  \label{errorerho}
\Eeq
for every $t\in[0,T]$.
Now, we deal with equations \eqref{discrmu} and~\eqref{prima}.
For the reader's convenience, by recalling that $\overg=g(\overrho)$ and $\underg=g(\underrho)$
(see~\eqref{defxingn}), we rewrite the former
in a different way, namely,
\Beq
  \bigl( 1+2g(\underrho) \bigr) \dt\hmu
  + \overmu \dt\hg - \Delta\overmu
  = 0 .
  \label{newdiscrmu}
\Eeq
Next, we take the difference between \eqref{newdiscrmu} and \eqref{prima}
and write it~as
\Bsist
  && \bigl( 1+2g(\hrho) \bigr) \, \dt(\hmu-\mu)
  - \Delta(\hmu-\mu)
  + (\hmu-\mu)
  \non
  \\
  && =
  - 2 \dt\mu \bigl( g(\hrho)-g(\rho) \bigr)
  - \dt\hg \, (\overmu-\mu)
  - \mu \, \dt \bigl( \hg-g(\rho) \bigr)
  \non
  \\
  && \quad {}
  + 2 \bigl( g(\hrho)-g(\underrho) \bigr) \dt\hmu
  -\Delta(\hmu-\overmu)
  + (\hmu-\mu) \,.
  \non
\Esist
Finally, we multiply {\juerg this} equality by~$(\hmu-\mu)$ and obtain the {\paolo following} identity:
\Bsist
  && \dt \bigl\{
    \bigl( \textstyle{\frac 12} + g(\hrho) \bigr) (\hmu-\mu)^2
  \bigr\}
  - \Delta(\hmu-\mu) \, (\hmu-\mu)
  + (\hmu-\mu)^2
  \non
  \\
  && = \dt g(\hrho) \, (\hmu-\mu)^2
  - 2 \dt\mu \, \bigl( g(\hrho)-g(\rho) \bigr) (\hmu-\mu)
  \non
  \\
  && \quad {}
  - \dt\hg \, (\overmu-\mu) (\hmu-\mu)
  - \mu \, \dt \bigl( \hg-g(\rho) \bigr) (\hmu-\mu)
  \non
  \\
  && \quad {}
  + 2 \bigl( g(\hrho)-g(\underrho) \bigr) \dt\hmu  \, (\hmu-\mu)
  -\Delta(\hmu-\overmu) \, (\hmu-\mu)
  + (\hmu-\mu)^2 .
  \qquad
  \label{testdiffmu}
\Esist
At this point, we integrate over~$Q_t$.
As $g$ is nonnegative, we~get:
\Beq
  \frac 12 \iO |(\hmu-\mu)(t)|^2 
  + \iot \normaV{(\hmu-\mu)(s)}^2 \, ds
  \leq \somma j17 I_j(t)\,,
  \label{inttestdiffmu}
\Eeq
with an obvious meaning of~$I_j(t)$, $j=1, \ldots,7$.
Now, we estimate {\juerg these} integrals, but the last one.
By combining the \holder, Young, and Sobolev, inequalities, and
in view of \eqref{stimaaux}, we~have that
\Bsist
  && I_1(t)
  \leq c \iot {\gianni\norma{\dt\hrho(s)}_4 \norma{(\hmu-\mu)(s)}_{{\ultime H}} \norma{(\hmu-\mu)(s)}_4} \, ds
  \non
  \\
  &&\leq c \iot {\gianni\norma{\dt\hrho(s)}_V \norma{(\hmu-\mu)(s)}_H \norma{(\hmu-\mu)(s)}_V} \, ds
  \non
  \\
   && {} \leq \sigma \iot \normaV{(\hmu-\mu)(s)}^2 \, ds
  + c_\sigma \iot \normaH{(\hmu-\mu)(s)}^2 \, ds\,, 
  \non
\Esist
where $\sigma>0$ is arbitrary. Similarly, we infer that
\Bsist
  && I_2(t) 
  \leq 2 \iot \norma{\dt\mu(s)}_{{\ultime H}} \norma{(\hrho-\rho)(s)}_4 \norma{(\hmu-\mu)(s)}_4 \, ds
  \non
  \\
  && \leq \sigma \iot \normaV{(\hmu-\mu)(s)}^2 \, ds
  + c_\sigma \iot \normaH{\dt\mu(s)}^2 \normaV{(\hrho-\rho)(s)}^2 \, ds .
  \non
\Esist
Notice that{\paolo , by means of the Gronwall lemma,} we shall be able to control the last integral 
 in terms 
of the $L^1(0,T)$-norm of the function $s\mapsto\normaH{\dt\mu(s)}^2$ 
(cf.~\eqref{regmu}). 
We use a similar procedure for the next integral 
and notice that the same remark holds, due to~\eqref{daauxbis}.
{\paolo Indeed, w}e have that
\Bsist
  && I_3(t)
  \leq \iot \norma{\dt\hg(s)}_\infty \normaH{(\hmu-\mu)(s)} \normaH{(\overmu-\mu)(s)} \, ds
  \non
  \\
  && \leq \iot \norma{\dt\hg(s)}_\infty \normaH{(\hmu-\mu)(s)} \bigl( \normaH{(\hmu-\mu)(s)} + \normaH{(\overmu-\hmu)(s)} \bigr) \, ds
  \non
  \\
  && \leq \iot \norma{\dt\hg(s)}_\infty^2 \normaH{(\hmu-\mu)(s)}^2 \, ds
  + c \iot \bigl( \normaH{(\hmu-\mu)(s)}^2 + \normaH{(\overmu-\hmu)(s)}^2 \bigr) \, ds
  \non
  \\
  && \leq c \iot \bigl( \norma{\dt\hg(s)}_\infty^2 + 1 \bigr) \normaH{(\hmu-\mu)(s)}^2 \, ds
  + c \tau^2\,,
  \non
\Esist
{\juerg where} the last inequality is due to \eqref{stimadiffmu}.
In order to treat~$I_4(t)$, we prove a preliminary estimate{\paolo ,
namely, that} 
\Bsist
  && |\dt \bigl( \hg-g(\rho) \bigr)| 
  \leq c \, \bigl\{
    |\overrho-\hrho| + |\underrho-\hrho| + |\hrho-\rho|
  \bigr\} \, |\dt\hrho|
  + c |\dt(\hrho-\rho)| 
  \qquad\qquad
  \label{perIquattro}
\Esist
{\paolo  \aeQ.} As we argue pointwise, we fix $(x,t)$ \aeQ\
and choose $n$ such that $t$ belongs to the interval $(n\tau,(n+1)\tau]$;
in order to simplify the notation,
we {\paolo omit writing at what point $(x,t)$ we work}.
By the mean value theorem, we find $r$ between $\rhon$ and $\rhonp$ such~that
\Bsist
  && \dt \bigl( \hg-g(\rho) \bigr)  
  = \frac {g(\rhonp)-g(\rhon)} \tau
  - g'(\rho) \dt\rho
  = g'(r) \, \frac {\rhonp-\rhon} \tau
  - g'(\rho) \dt\rho
  \non
  \\
  && = g'(r) \dt\hrho
  - g'(\rho) \dt\rho
  = \bigl( g'(r) - g'(\rho) \bigr) \dt\hrho
  + g'(\rho) \bigl( \dt\hrho - \dt\rho \bigr) .
  \non
\Esist
As $g'$ is bounded and Lipschitz continuous, we infer~that
\Beq
  |\dt \bigl( \hg-g(\rho) \bigr)| 
  \leq c |r-\rho| \, |\dt\hrho|
  + c |\dt\hrho - \dt\rho| .
  \non
\Eeq
{\gianni
On the other hand, we have
\Bsist
  && |r-\rho|
  \leq |r-\rhon| + |\rhon-\rho|
  \leq |\rhonp-\rhon| + |\rhon-\rho|
  \non
  \\
  && = |\overrho-\underrho| + |\underrho-\rho|
  \leq |\overrho-\hrho| + 2 |\underrho-\hrho| + |\hrho-\rho| . 
  \non
\Esist
Hence, \eqref{perIquattro} follows, and we can use it to estimate~$I_4(t)$.
We also account for the boundedness of~$\mu$ and for identity~\eqref{interpLinfty}
and the analogue {\paolo identity} concerning~$\underz$.
We~have:
\Bsist
  && I_4(t)
  \leq c \intQt \bigl\{
    |\overrho-\hrho| + |\underrho-\hrho| + |\hrho-\rho|
  \bigr\} \, |\dt\hrho| \, |\hmu-\mu|
  \non
  \\  
  && \quad {}
  + c \intQt |\dt(\hrho-\rho)| \, |\hmu-\mu|
  \non
  \\
  && \leq c \iot \bigl\{
    \normaH{\overrho(s)-\hrho(s)}^2 + \normaH{\underrho(s)-\hrho(s)}^2 + \normaH{\hrho(s)-\rho(s)}^2
  \bigr\} \, \norma{\dt\hrho(s)}_\infty^2 \, ds
  \non
  \\
  && \quad {}
  + \frac12 \iot \normaH{\dt(\hrho-\rho)(s)}^2 \, ds
  + c \iot \normaH{(\hmu-\mu)(s)}^2 \, ds 
  \non
  \\
  && \leq c \, \tau^2 \norma{\dt\hrho}_{\L\infty H}^2 \, \norma{\dt\hrho}_{\L2\Linfty}^2 
  + c \iot \norma{\dt\hrho(s)}_\infty^2 \, \normaH{\hrho(s)-\rho(s)}^2 \, ds
  \non
  \\
  && \quad {}
 + \frac12 \iot \normaH{\dt(\hrho-\rho)(s)}^2 \, ds
  + c \iot \normaH{(\hmu-\mu)(s)}^2 \, ds \,;
  \non
\Esist
{\paolo furthermore,}  estimates \eqref{stimedthat} and \eqref{daaux} allow us to infer~that
\Bsist
  && I_4(t)
  \leq c \, \tau^2
  + c \iot \norma{\dt\hrho(s)}_\infty^2 \, \normaV{\hrho(s)-\rho(s)}^2 \, ds
  + \frac12 \iot \normaH{\dt(\hrho-\rho)(s)}^2 \, ds
  \non
  \\
  && \quad {}
  + c \iot \bigl( 1 + \norma{\dt\hrho(s)}_\infty^2 \bigr) \normaH{(\hmu-\mu)(s)}^2 \, ds .
  \non
\Esist
}%
Next, we deal with~$I_5(t)$.
By accounting for~\eqref{interpLdue} and~\eqref{daaux}, we deduce~that
\Bsist
  && I_5(t)
  \leq c \iot \norma{(\hrho-\underrho)(s)}_\infty \norma{\dt\hmu(s)}_{{\ultime H}} \norma{(\hmu-\mu)(s)}_{{\ultime H}} \, ds
  \non
  \\
  && \leq c \, \tau^2 \norma{\dt\hrho}_{\L2\Linfty}^2
  + c \iot \normaH{\dt\hmu(s)}^2 \normaH{(\hmu-\mu)(s)}^2 \, ds
  \non
  \\
  && \leq c \, \tau^2 
  + c \iot \normaH{\dt\hmu(s)}^2 \normaH{(\hmu-\mu)(s)}^2 \, ds .
  \non
\Esist
We note at once that we shall be able to control 
even the last terms of the last two estimates
with the help of the Gronwall lemma,
in view of~\eqref{daaux} and~\eqref{stimedthat}, respectively.
Finally, thanks to~\eqref{interpLdue} once more, we~have:
\Bsist
  && I_6(t)
  = \intQt \nabla(\hmu-\overmu) \cdot \nabla(\hmu-\mu) 
  \non
  \\
  && \leq \sigma \intQt |\nabla(\hmu-\mu)|^2
  + c_\sigma \intQt |\nabla(\hmu-\overmu)|^2 
  \non
  \\
  \noalign{\allowbreak}
  && \leq \sigma \intQt |\nabla(\hmu-\mu)|^2
  + c_\sigma \, \tau^2 \norma{\nabla\dt\hmu}_{\L2H}^2  
  \non
  \\
  && \leq \sigma \intQt |\nabla(\hmu-\mu)|^2
  + c_\sigma \, \tau \,,
  \label{stimaIsei}
\Esist
{\juerg where} the last inequality is a consequence of 
the non-sharp estimate~\eqref{nonsharp}.
We stress that $I_6$ is the only term
of order $\tau$ instead of~$\tau^2$.
At this point, we collect all the estimates 
of the integrals~$I_j$ we have obtained,
and come back to \accorpa{testdiffmu}{inttestdiffmu}.
If we choose $\sigma$ small enough,  
we conclude~that
{\gianni
\Bsist
  && \frac12 \,\normaH{(\hmu-\mu)(t)}^2
  + \frac12 \iot \normaV{(\hmu-\mu)(s)}^2 \, ds
  \non
  \\
  && {}\leq c \, 
  \Bigl\{
    \tau
    + \iot  \bigl( 1+ \norma{\dt\hg(s)}_\infty^2 + \normaH{\dt\hmu(s)}^2  \bigr)
                   \normaH{(\hmu-\mu)(s)}^2 \, ds 
   \non
   \\
    && \qquad {}
   + \iot \bigl( \normaH{\dt\mu(s)}^2 + \norma{\dt\hrho(s)}_\infty^2 \bigr) \, 
   \normaV{\hrho(s)-\rho(s)}^2 \, ds
  \non
  \\
  && \qquad {}
    + \intQt |\dt(\hrho-\rho)|^2 
  \Bigr\} + \frac12 \iot \normaH{\dt(\hrho-\rho)(s)}^2 \, ds
  \qquad
  \label{erroremu}
\Esist
}%
for every $t\in[0,T]$. 
Now, we {\paolo revert} to~\eqref{errorerho} and add it to~\eqref{erroremu}.
After rearranging, 
we apply the Gronwall lemma and obtain~\eqref{errore}.
This concludes the proof.


\section{Proof of Theorem \ref{Errorebis}}
\label{ProofErrorebis}
\setcounter{equation}{0}

As  is clear from the proof of Theorem~\ref{Errore}{\paolo , 
 to obtain estimate~\eqref{errorebis} there {\pier is}
just one step to modify,} namely,
the estimate of~$I_6$ (see~\eqref{stimaIsei}), 
which was based on the non-sharp inequality~\eqref{nonsharp}.
Thus, we only have to prove that our further assumption \eqref{hperrorebis}
implies {\paolo that }$I_6$ {\paolo must} be of order $\tau^2$, {\paolo not}~$\tau$.
Moreover, it is clear that this is true whenever we improve~\eqref{nonsharp}
and replace it~by
\Beq
  \norma{\nabla\dt\hmu}_{\L2H}^2 \leq c ,
  \quad \hbox{i.e.,} \quad
  \tau \somma n0{N-1} \normaH{\nabla\dtau\mun}^2 \leq c .
  \label{sharp}
\Eeq
Hence, it suffices to prove~\eqref{sharp}.
In order to make our argument transparent, 
we {\juerg prove some additional} estimates{\paolo , the first of which} holds under assumption~\eqref{hperrorebis}.

\step
Further a priori estimates

We prepare an estimate of $\normaH{\dtau\muz}$.
In view of~\eqref{hperrorebis},
we write equation \eqref{dprima}, with $n=0$, in the form:
\Beq
  (1+2\gz) \dtau\muz - \tau \Delta\dtau\muz
  = \Delta\muz - \mu_1 \, \dtau\gz,
  \non
\Eeq
and test {\paolo it} by $\dtau\muz$.
As $\gz$ is nonnegative, we immediately {\paolo arrive at}
\Beq
  \iO |\dtau\muz|^2
  + \tau \iO |\nabla\dtau\muz|^2
  \leq \bigl( \norma{\Delta\muz}_{{\ultime H}} + \norma{\mu_1}_4 \, \norma{\dtau\rhoz}_4 \bigr) \, \norma{\dtau\muz}_{{\ultime H}} \,.
  \non
\Eeq
Thanks to~\eqref{hperrorebis}, the Sobolev inequality,
\eqref{sestastimad}, and \eqref{stimaauxd},
we~deduce~that
\Beq
  \normaH{\dtau\muz}
  + \tau \iO |\nabla\dtau\muz|^2
  \leq c .
  \label{stimadtmuz}
\Eeq

Let us come to the basic estimate we need.
We improve \eqref{quartastimad} and obtain a bound
for the second{\paolo -}difference quotients~$\dtau^2\rhon$ (see~\eqref{defdeltadue}).
We write \eqref{dseconda}, with $(n+1)$ in place of $n$,
and test the difference between the {\juerg resulting} {\paolo relation} and \eqref{dseconda} itself
by $(\dtau\rhonp-\dtau\rhon)$.
We~find:
\Bsist
  && \iO |\dtau\rhonp-\dtau\rhon|^2
  + \tau \iO \nabla\dtau{\gianni\rhonp} \cdot \nabla(\dtau\rhonp-\dtau\rhon)
  \non
  \\
  && {\gianni = - \iO \bigl( f'(\rhonpp) - f'(\rhonp) \bigr) (\dtau\rhonp-\dtau\rhon) }
  \non
  \\
  && \quad {}
  + \tau \iO \bigl( g'(\rhonp)\dtau\mun+\mun\dtau(g'(\rhon)) \bigr) (\dtau\rhonp-\dtau\rhon)
  \non
  \\
  && \leq C \, \tau \iO \bigl( |\dtau\rhonp| + |\dtau\mun| + |\mun| \, |\dtau\rhon| \bigr) |\dtau\rhonp-\dtau\rhon| \,.
  \label{perauxdue}
\Esist
By the elementary identity \eqref{elementare}, we have:
\Beq
  \nabla\dtau{\gianni\rhonp} \cdot \nabla(\dtau\rhonp-\dtau\rhon)
  = \frac 12 |\nabla\dtau\rhonp|^2
  - \frac 12 |\nabla\dtau\rhon|^2
  + \frac 12 |\nabla\dtau\rhonp-\nabla\dtau\rhon|^2 .
  \non
\Eeq
On the other hand, by the Sobolev inequality, \eqref{sestastimad}, and 
\eqref{stimaauxd}, we infer that
\Bsist
  && C \, \tau \iO \bigl( {\gianni|\dtau\rhonp|} + |\dtau\mun| + |\mun| \, |\dtau\rhon| \bigr) |\dtau\rhonp-\dtau\rhon|
  \non
  \\
  && \leq \frac 1 2 \iO |\dtau\rhonp-\dtau\rhon|^2
  + c \, \tau^2 \bigl( {\gianni\norma{\dtau\rhonp}_{{\ultime H}}^2} + \norma{\dtau\mun}_{{\ultime H}}^2 + \norma\mun_4^2 \, \norma{\dtau\rhon}_4^2 \bigr)
  \non
  \\
  && \leq \frac 1 2 \iO |\dtau\rhonp-\dtau\rhon|^2
  + c \, \tau^2 \bigl( {\gianni\norma{\dtau\rhonp}_{{\ultime H}}^2} + \normaH{\dtau\mun}^2 + \normaV\mun^2 \, \normaV{\dtau\rhon}^2 \bigr) 
  \non
  \\
  && \leq \frac 1 2 \iO |\dtau\rhonp-\dtau\rhon|^2
  + c \, \tau^2 \bigl( {\gianni\normaH{\dtau\rhonp}^2} + \normaH{\dtau\mun}^2 + 1 \bigr) .
  \non
\Esist
Now, we {\paolo combine} {\juerg this} estimate, the identity {\paolo just above}, and \eqref{perauxdue}.
Then, we divide by~$\tau$ and sum over $n=0,\dots,m-1$, where $1\leq m\leq N-1$.
We obtain:
\Bsist
  && \frac \tau 2 \, \somma n0{m-1} \normaH{\dtau^2\rhon}^2
  + \frac 12 \, \normaH{\nabla\dtau\rhom}^2
  + \frac 12 \, \somma n0{m-1} \normaH{\nabla\dtau\rhonp-\nabla\dtau\rhon}^2
  \non
  \\
  && \leq \frac 12 \, \normaH{\nabla\dtau\rhoz}^2
  + {\gianni c \, \tau \somma n0{m-1} \normaH{\dtau\rhonp}^2}
  + c \, \tau \somma n0{m-1} \normaH{\dtau\mun}^2
  + c .
  \non
\Esist
At this point, by \eqref{stimaauxd}, \eqref{secondastimad}, 
and \eqref{sestastimad}, we conclude~that
\Beq
  \tau \somma n0{N-2} \normaH{\dtau^2\rhon}^2
  \leq c .
  \label{auxdue}
\Eeq

\step
Consequence

{\paolo With a view toward deriving an} estimate for $\dtau^2\gn$, {\paolo we begin by arguing} pointwise. So, for a.a.\ $(x,t)\in Q$
({\paolo once again we omit writing at what point of $Q$ we work}) and for suitable $r_1$ between $\rhonpp$ and~$\rhonp$,
and $r_2$ between $\rhon$ and~$\rhonp$, we~have by the Taylor formula: 
\Bsist
  && |\dtau^2\gn|
  = \tau^{-2} \, |g(\rhonpp) - g(\rhonp) + g(\rhon) - g(\rhonp)|
  \non
  \\
  && = \tau^{-2} \, | g'(\rhonp) (\rhonpp-\rhonp)
  + \textstyle{\frac 12} \, g''(r_1) (\rhonpp-\rhonp)^2 
  \non
  \\
  && \quad {}
  + g'(\rhonp) (\rhon-\rhonp)
  + \textstyle{\frac 12} \, g''(r_2) (\rhon-\rhonp)^2 |
  \non
  \\
  && \leq c \, |\dtau^2\rhon|
  + c \, \bigl( |\dtau\rhonp|^2 + |\dtau\rhon|^2 \bigr) .
  \non
\Esist
Now, we square {\juerg this} pointwise estimate, integrate over~$\Omega$, 
sum over~$n$, and deduce {\juerg that}
\Beq
  \tau \somma n0{N-2} \normaH{\dtau^2\gn}^2
  \leq c \, \tau \somma n0{N-2} \normaH{\dtau^2\rhon}^2
  + c \, \tau \somma n0{N-1} \norma{\dtau\rhon}_4^4 \,.
  \non
\Eeq
Then, \eqref{auxdue}, the Sobolev inequality, and \eqref{stimaauxd}  yield:
\Beq
  \tau \somma n0{N-2} \normaH{\dtau^2\gn}^2 \leq c .
  \label{daauxdue}
\Eeq

\step
Proof of Theorem~\ref{Errorebis}

As said before, it suffices to prove~\eqref{sharp}.
{\paolo We reason that,} in order to obtain the analogous estimate for the solution to the continuous problem,
one first differentiates \eqref{prima} with respect to time
and then tests the resulting equality by~$\dt\mu$;
this yields the desired term $\displaystyle{\intQt|\nabla\dt\mu|^2}$ on the \lhs. 
{\paolo The idea is} to perform the corresponding procedure on the discrete equation~\eqref{dprima}.
However, it turns out that the calculation in the discrete case becomes simpler
if one tests by the analogue of the product $(1+2g(\rho))\dt\mu$.
To simplify the notation,
we introduce the vector $\pi$ defined by 
\Beq
  \pin := (1+2\gn) \dtau\mun
  \quad \hbox{for $n=0,\dots,N-1$}.
  \label{defpin}
\Eeq
We write \eqref{dprima} with $(n+1)$ in place of~$n$,
and take the difference between the resulting equality and \eqref{dprima} itself.
Then, we test {\juerg this} difference by $\pinp$ and integrate over~$\Omega$.
By {\paolo taking} the elementary identity~\eqref{elementare} {\paolo into account},
we obtain for $n=0,\dots,N-2$ {\juerg that}
\Bsist
  && \frac 12 \iO |\pinp|^2
  - \frac 12 \iO |\pin|^2
  + \frac 12 \iO |\pinp-\pin|^2
  + \tau \iO \nabla\dtau\munp \cdot \nabla\pinp
  \non
  \\
  && = - \iO \bigl( \munpp \, \dtau\gnp - \munp \, \dtau\gn \bigr) \, \pinp  \,.
  \non
\Esist
By computing the {\juerg fourth} term on the \lhs\ with the help of~\eqref{defpin}, 
recalling that $g$ is nonnegative, and rearranging, we deduce that
\Bsist
  && \iO |\pinp|^2
  - \iO |\pin|^2
  + \iO |\pinp-\pin|^2
  + 2\tau \iO |\nabla\dtau\munp|^2
  \non
  \\
  \noalign{\allowbreak}
  && \leq - 2 \iO \bigl( \munpp \, \dtau\gnp - \munp \, \dtau\gn \bigr) \, \pinp 
  - 4\tau \iO \dtau\munp \, \nabla\dtau\munp \cdot \nabla\gnp 
  \qquad
  \non
  \\
  && = - 2 \iO \pinp (\munpp-\munp) \, \dtau\gnp
  - 2 \iO \pinp \, \munp \bigl( \dtau\gnp-\dtau\gn \bigr)
  \non
  \\
  && \quad {}
  - 4\tau \iO \dtau\munp \, \nabla\dtau\munp \cdot \nabla\gnp
  \non
  \\
  && = - 2\tau \iO \pinp \, \dtau\munp \, \dtau\gnp
  - 2\tau \iO \pinp \, \munp \, \dtau^2 \gn
  \non
  \\
  && \quad {}
  - 4\tau \iO \dtau\munp \, \nabla\dtau\munp \cdot \nabla\gnp \,.
  \label{perstimabis}
\Esist
Now, we estimate each term of the \rhs\ separately, before summing over~$n$,
in order to simplify the notation.
For the first one,
we use \holder\ and Sobolev inequalities, and estimates~\eqref{stimaauxd} and~\eqref{auxdue}.
We~have:
\Bsist
  && - 2\tau \iO \pinp \, \dtau\munp \, \dtau\gnp
  \non
  \\
  && \leq 2\tau \norma\pinp_4 \, \norma{\dtau\munp}_{{\ultime H}} \, \norma{\dtau\gnp}_4 
  \leq c \, \tau \normaV{\dtau\munp} \, \normaH{\dtau\munp} \, \normaV{\dtau\rhonp}
  \non
  \\
  && \leq \frac \tau 4 \, \bigl( \normaH{\nabla\dtau\munp}^2 + \normaH{\dtau\munp}^2 \bigr)
  + c \, \tau \normaH{\dtau\munp}^2 
  \non
  \\
  && \leq \frac \tau 4 \, \normaH{\nabla\dtau\munp}^2 
  + c \, \tau  \, \normaH\pinp^2  
  \non
  \\
  && \leq \frac \tau 4 \, \normaH{\nabla\dtau\munp}^2 
  + c \, \tau  \, \normaH\pin^2
  + C \, \tau  \, \normaH{\pinp-\pin}^2\,.
  \non
\Esist
{\paolo For} $\tau$ small enough, namely, {\paolo for} $\tau\leq 1/(3C)$, we conclude that
\Beq
  - 2\tau \iO \pinp \, \dtau\munp \, \dtau\gnp
  \leq \frac \tau 4 \, \normaH{\nabla\dtau\munp}^2 
  + c \, \tau  \, \normaH\pin^2
  + \frac 13 \, \normaH{\pinp-\pin}^2 .
  \non
\Eeq
Next, by~\eqref{sestastimad}, we similarly have:
\Bsist
  && - 2\tau \iO \pinp \, \munp \, \dtau^2 \gn
  \leq c \, \tau \norma{\dtau\munp}_4 \, \norma\munp_4 \, \norma{\dtau^2\gn}_{{\ultime H}}
  \non
  \\
  && \leq \frac \tau 4 \, \bigl( \normaH{\nabla\dtau\munp}^2 + \normaH{\dtau\munp}^2 \bigr)
  + c \, \tau \, \normaV\munp^2 \, \normaH{\dtau^2\gn}^2
  \non
  \\
  && \leq \frac \tau 4 \, \normaH{\nabla\dtau\munp}^2
  + \tau \, \normaH{\pin}^2
  + \frac 13 \, \normaH{\pinp-\pin}^2 
  + c \, \tau \, \normaH{\dtau^2\gn}^2 \,,
  \non
\Esist
{\paolo {\pier for} sufficiently small $\tau$}. Finally, by accounting for \eqref{quintastimad},
Sobolev inequality, and the compactness inequality~\eqref{compact},
we have that, for $\tau$ small enough,
\Bsist
  && - 4\tau \iO \dtau\munp \, \nabla\dtau\munp \cdot \nabla\gnp
  \leq c \, \tau \norma{\dtau\munp}_4 \, \norma{\nabla\dtau\munp}_{{\ultime H}} \, \norma{\nabla\rhonp}_4
  \non
  \\  
  && \leq \frac \tau 4 \, \norma{\nabla\dtau\munp}_{{\ultime H}}^2 
  + c \, \tau \norma{\dtau\munp}_4^2 \, \normaV{\nabla\rhonp}^2
  \leq \frac \tau 4 \, \norma{\nabla\dtau\munp}_{{\ultime H}}^2 
  + c \, \tau \norma{\dtau\munp}_4^2 
  \non
  \\  
  && \leq \frac \tau 4 \, \norma{\nabla\dtau\munp}_{{\ultime H}}^2 
  + \tau \bigl( \textstyle{\frac 14} \, \normaH{\nabla\dtau\munp}^2 + c \normaH{\dtau\munp}^2 \bigr)
  \non
  \\
  && \leq \frac \tau 2 \, \norma{\nabla\dtau\munp}_{{\ultime H}}^2 
  + c \, \tau \normaH\pinp^2
  \leq \frac \tau 2 \, \norma{\nabla\dtau\munp}_{{\ultime H}}^2 
  + c \, \tau \normaH\pin^2
  + \frac 13 \, \normaH{\pinp-\pin}^2 \,.
  \non
\Esist
At this point, we combine the inequalities just obtained with \eqref{perstimabis}
and note that the terms involving $\pinp-\pin$ cancel out.
Then, we sum over $n=0,\dots,(m-1)$, with $1\leq m\leq (N-1)$.
We obtain:
\Beq
  \normaH\pim^2 
  + \tau \somma n0{m-1} \normaH{\nabla\dtau\munp}^2
  \leq \normaH\piz^2
  + \tau \somma n0{m-1} \normaH\pin^2
  + \tau \somma n0{m-1} \normaH{\dtau^2\gn}^2
  \non
\Eeq
and the discrete Gronwall lemma allows us to deduce that 
\Beq
  \normaH\pim^2 
  + \tau \somma n0{m-1} \normaH{\nabla\dtau\munp}^2
  \leq c \Bigl(
    \normaH\piz^2
    + c \, \tau \somma n0{N-2} \normaH{\dtau^2\gn}^2    
  \Bigr) 
  \non
\Eeq
for $1\leq m\leq (N-1)$.
From \eqref{stimadtmuz} and \eqref{daauxdue},
we infer that
\Beq
  \tau \somma n0{N-2} \normaH{\nabla\dtau\munp}^2
  \leq c .
  \non
\Eeq
This and \eqref{stimadtmuz} yield \eqref{sharp}, and the proof is complete.

\Brem
\label{Piureg}
As a consequence of  estimates~\eqref{sharp} and~\eqref{auxdue},
the solution to the continuous problem enjoys the {\paolo following additional} regularity properties:
\Beq
  \nabla\dt\mu \in \LQ2
  \aand
  \dt^2\rho \in \LQ2\,. 
  \non
\Eeq
This can give even more: for instance, equation~\eqref{seconda} can be differentiated with respect to time{\paolo , to show that} 
$\Delta\dt\rho$ belongs to $\LQ2$ as well, {\paolo so as to} conclude that
\Beq
  \rho \in \H2H \cap \H1W.
  \non
\Eeq
However, {\juerg this} regularity result could be proved {\paolo formally and directly for} the continuous problem.
\Erem



\bibliographystyle{abbrv}
\bibliography{CoGiKrPoSp}
%
\End{document}

\bye